\documentclass[12pt,leqno]{amsart}

\usepackage{amssymb}
\usepackage{amsthm}
\usepackage{amsmath}
\usepackage{amscd}
\usepackage[mathscr]{eucal}
\usepackage{verbatim}

\numberwithin{equation}{section}

\theoremstyle{plain}
\newtheorem{theorem}{Theorem}[section]
\newtheorem{corollary}[theorem]{Corollary}
\newtheorem{lemma}[theorem]{Lemma}
\newtheorem{proposition}[theorem]{Proposition}

\theoremstyle{definition}
\newtheorem{definition}[theorem]{Definition}
\newtheorem{remark}[theorem]{Remark}

\theoremstyle{remark}
\newtheorem*{notation}{Notation}

\newcommand{\R}{\mathbb{R}}
\newcommand{\Q}{\mathbb{Q}}
\newcommand{\Z}{\mathbb{Z}}
\newcommand{\N}{\mathbb{N}}
\newcommand{\C}{\mathbb{C}}
\newcommand{\h}{\mathbb{H}}
\renewcommand{\H}{\mathbb{H}}

\newcommand{\G}{\Gamma}
\newcommand{\g}{\gamma}

\newcommand{\la}{\lambda}

\newcommand{\x}{\mathbf{x}}

\newcommand{\back}{\backslash}

\newcommand{\V}{V(\Q)}

\newcommand{\zxz}[4]{\begin{pmatrix} #1 & #2 \\ #3 & #4 \end{pmatrix}}

\newcommand{\kzxz}[4]{\left(\begin{smallmatrix} #1 & #2 \\ #3 & #4\end{smallmatrix}\right) }
\newcommand{\kabcd}{\kzxz{a}{b}{c}{d}}

\newcommand{\re}{\operatorname{Re}}
\newcommand{\Res}{\operatorname{Res}}

\newcommand{\calE}{\mathcal{E}}
\newcommand{\calF}{\mathcal{F}}

\newcommand{\calM}{\mathcal{M}}

\newcommand{\calQ}{\mathcal{Q}}

\newcommand{\calZ}{\mathcal{Z}}

\newcommand{\frake}{\mathfrak e}

\newcommand{\eps}{\varepsilon}
\newcommand{\bs}{\backslash}

\newcommand{\vol}{\operatorname{vol}}
\newcommand{\tr}{\operatorname{tr}}

\newcommand{\Span}{\operatorname{span}}

\newcommand{\Sl}{\operatorname{SL}}

\newcommand{\SL}{\operatorname{SL}}

\newcommand{\Spin}{\operatorname{Spin}}
\newcommand{\PSL}{\operatorname{PSL}}
\newcommand{\Symp}{\operatorname{Sp}}
\newcommand{\Mp}{\operatorname{Mp}}
\newcommand{\Orth}{\operatorname{O}}

\newcommand{\Aut}{\operatorname{Aut}}

\newcommand{\SO}{\operatorname{SO}}

\newcommand{\End}{\operatorname{End}}

\newcommand{\Ei}{\operatorname{Ei}}
\newcommand{\erfc}{\operatorname{erfc}}
\newcommand{\Iso}{\operatorname{Iso}}

\newcommand{\cha}{\operatorname{\widehat{CH}}}

\begin{document}

\title[Traces of CM values of modular functions]{Traces of CM values of modular functions}

\author[Jan H.~Bruinier and Jens Funke]{Jan Hendrik Bruinier and Jens Funke*}

\address{Mathematisches Institut, Universit\"at zu K\"oln, Weyertal 86--90, D-50931 K\"oln, Germany}
\email{bruinier@math.uni-koeln.de } 
\address{Department of Mathematical Sciences, New Mexico State University, P.O.~Box 30001, 3MB, Las Cruces, NM 88001, USA}
\email{jfunke@nmsu.edu}

\thanks{* Partially supported by NSF-grant DMS-0305448 and NSF-grant DMS-0211133 at the Fields Institute, Toronto}


\date{July 29, 2004}

\begin{abstract}
Zagier proved that the traces of singular moduli, i.e., the sums of the values of the classical $j$-invariant over quadratic irrationalities, are the Fourier coefficients of a modular form of weight $3/2$ with poles at the cusps. 
Using the theta correspondence, we generalize this result to traces of CM values of (weakly holomorphic) modular functions on modular curves of arbitrary genus. 
We also study the theta lift for the weight $0$ Eisenstein series for $\Sl_2(\Z)$ and realize a certain generating series of arithmetic intersection numbers as the 
derivative of Zagier's Eisenstein series of weight $3/2$. This recovers a result of Kudla, Rapoport and Yang. 
\end{abstract}

\maketitle

\section{Introduction}


In \cite{ZagierTr}, Zagier considers the normalized Hauptmodul $J(z) = j(z) -744$ for the group $\Gamma(1)=\PSL_2(\Z)$, where $j(z) = e^{- 2\pi i z} +744 + 196884 e^{2\pi i z} + \dots $ is the classical $j$-invariant on the complex upper half plane $\h$. 
Let $D$ be a positive integer and write $\calQ_{D}$ for the set of positive definite integral binary quadratic forms $[a,b,c]$ of discriminant $-D=b^2-4ac$.
The group $\Gamma(1)$ acts on $\calQ_{D}$. If $Q=[a,b,c]\in \calQ_{D}$ we write $\Gamma(1)_Q$ for the stabilizer of $Q$ in $\Gamma(1)$ and  
$\alpha_Q=\frac{-b+i\sqrt{D}}{2a}$ for the corresponding CM point in $\H$.
The values of $j$ at such points $\alpha_Q$ are known as singular moduli. They play an important role in many branches of number theory.
The modular trace of $J$ of index $D$ is defined as
\begin{equation}\label{modtrace}
{\bf t}_J(D) = \sum_{Q\in \calQ_{D}/\Gamma(1)} \frac{1}{|\Gamma(1)_Q|} 
J(\alpha_Q).
\end{equation}
By the theory of complex multiplication, ${\bf t}_J(D)$ can also be viewed as a suitable Galois trace. It is a rational integer.

Zagier shows that the generating series
\begin{equation}\label{Zagier1}
-q^{-1} +2 + \sum_{D=1}^{\infty} {\bf t}_J(D) q^D= -q^{-1} +2 -248 q^3+492 q^4-4119q^7+7256q^8+\dots
\end{equation}
is a meromorphic modular form of weight $3/2$ for the Hecke subgroup $\G_0(4)$ whose poles are supported at the cusps. Here $q = e^{2\pi i\tau}$ with $\tau = u+iv \in \h$. He gives two proofs of this result. The first uses certain recursion relations for the ${\bf t}_J(D)$, the second uses Borcherds products on $\Sl_2(\Z)$ and an application of Serre duality. Both proofs 
rely on the fact that (the compactification of) $\G(1) \back \h$ has genus zero. 
In  \cite{Kim1,Kim2}, Kim extended Zagier's results to other modular curves of genus zero using similar methods. 

It is quite interesting to compare this result with an older theorem of Zagier \cite{Zagier} concerning the Hurwitz-Kronecker class numbers $H(D) = \sum_{Q\in \calQ_{D}/\Gamma(1)} \frac{1}{|\Gamma(1)_Q|}$, which we consider here as the trace ${\bf t}_1(D)$ of the constant modular function $1$ of weight $0$. Zagier constructs a certain Eisenstein series $\calF(\tau,s)$ of weight $3/2$ and shows that for the special value at $s = \tfrac12$ (in our normalization)  
\begin{equation}\label{Zagier2}
\calF(\tau,\frac12) = \sum_{D =0}^{\infty}  {\bf t}_1(D)   q^D + \frac{1}{16 \pi\sqrt{v}} \sum_{N=- \infty}^{\infty} \beta(4\pi N^2 v) q^{-N^2}
\end{equation}
is a non-holomorphic modular form of weight $3/2$ for $\G_0(4)$. Here ${\bf t}_1(0) = -\tfrac1{12}= \vol(\Gamma(1) \back  \h)$ and $\beta(s) = \int_1^{\infty} t^{-3/2} e^{-st} dt$. It is striking that while the positive Fourier coefficients of \eqref{Zagier1} and \eqref{Zagier2} are both traces of modular functions, the negative coefficients are very different in nature. Furthermore, Zagier's proofs for \eqref{Zagier1} and \eqref{Zagier2} are totally different.

In \cite{FCompo}, the second named author extended \eqref{Zagier2} to realize the generating series of the class numbers of CM points for general congruence subgroups $\G$ as the holomorphic part of a non-holomorphic modular form of weight $3/2$. These modular forms take the same form as in \eqref{Zagier2} and are obtained as a theta integral
\begin{equation}\label{thetaintro}
I(\tau,1) = \int_{\G \back \h} 1 \cdot \theta_{L}(\tau,z,\varphi) \,\frac{dx\,dy}{y^2},
\end{equation}
integrating the constant function $1$ against a theta series associated to an even lattice $L$ of signature $(1,2)$ and a certain Schwartz function $\varphi$ coming from \cite{KMI}. 

\smallskip

In the present paper, we use the method of \cite{FCompo} to generalize $\eqref{Zagier1}$ to traces ${\bf t}_f$ of arbitrary modular functions $f$ of weight $0$ whose poles are supported at the cusps on modular curves of higher genus. 
Namely, we consider the theta integral $I(\tau,f)$ replacing in \eqref{thetaintro} the constant $1$ by the more general modular function $f$. Here the starting point is that $I(\tau,f)$ does converge since the decay of the theta kernel turns out to be faster than the exponential growth of $f$ at the cusps, see Proposition~\ref{thetagrowth}. Furthermore, the Schwartz function $\varphi$ underlying the theta kernel is closely related to a Green function for the CM points constructed by Kudla \cite{KAnn}.
This approach also gives a unifying proof for \eqref{Zagier1} and \eqref{Zagier2}. Furthermore, we obtain geometric interpretations for the constant and the negative Fourier coefficients.
For instance, the constant coefficient can be interpreted as  the ``average value''
\[
-\frac{1}{2\pi}\int_{\G \back \h}^{reg} f(z) \,\frac{dx\,dy}{y^2},
\]
of $f$ on $\G \bs \H$.
Here $\int_{\G \back \h}^{reg}$ indicates a certain kind of regularization of the divergent integral. The negative coefficients involve data coming from infinite geodesics joining two cusps of $\G\back \h$.

\smallskip

To illustrate our result, we now describe a special case, see section~\ref{sec:ex}. For the general statement which is phrased in terms of the orthogonal group of a rational quadratic space of signature $(1,2)$, see Theorem~\ref{MAIN}.  

Let $p$ be a prime (or $p=1$). For a positive integer $D$, we consider
the subset $\calQ_{D,p}$ of quadratic forms $[a,b,c]\in \calQ_{D}$ such that $a\equiv 0\pmod{p}$. Note that $\G_0^{\ast}(p)$, the extension of the Hecke group $\Gamma_0(p)\subset\Gamma(1)$ with the Fricke involution 
 $W_p=\kzxz{0}{-1}{p}{0}$, acts on $\calQ_{D,p}$ with finitely many orbits.

Let $f$ be a modular function (of weight $0$) for $\Gamma^{\ast}_0(p)$ whose poles are supported at the cusp and denote its Fourier expansion by
$f(z)=\sum_{n\gg-\infty}a(n)e(nz)$.
We define the modular trace of $f$ of index $D$ by
\begin{equation}\label{modtrace2}
{\bf t}^{\ast}_f(D) = \sum_{Q\in \calQ_{D,p}/\Gamma^{\ast}_0(p)} \frac{1}{|\Gamma^{\ast}_0(p)_Q|} 
f(\alpha_Q),
\end{equation}
where $\Gamma^{\ast}_0(p)_Q$ is the stabilizer of $Q$ in $\Gamma^{\ast}_0(p)$. Finally, we put $\sigma_1(0)=-\frac{1}{24}$ and $\sigma_1(n)=\sum_{t\mid n} t$ for $n\in \Z_{\geq0}$ and $\sigma_1(x) =0$ for $x \notin \Z_{\geq0}$.

\begin{theorem} \label{thm:intro}
Let $f$ be a modular function for $\Gamma^{\ast}_0(p)$ and denote its Fourier expansion as above.
Assume that the constant coefficient $a(0)$ vanishes.
Then
\begin{align*}
G(\tau,f)&=\sum_{\substack{D>0}}  {\bf t}^{\ast}_f(D) q^D +
\sum_{n\geq 0}\big(\sigma_1(n)  + p \sigma_1(n/p) \big) a(-n) \\
& \phantom{=}{}-\sum_{m>0}\sum_{n>0}m a(-mn)q^{-m^2}
\end{align*}
is a meromorphic modular form of weight $3/2$, holomorphic outside the cusps, for the group $\G_0(4p)$ satisfying the Kohnen plus space condition (see \eqref{kohnencond}).
If $a(0)$ does not vanish, then in addition non-holomorphic terms as in \eqref{Zagier2} occur.

\end{theorem}
 
For $p=1$, and  $f=J$, we recover \eqref{Zagier1}. 



\medskip

One can also consider the theta lift $I(\tau,f)$ for other types of automorphic forms of weight $0$. We consider $I(\tau,\calE_0(z,s))$, where $\calE_0(z,s)$ is the (normalized) Eisenstein series for $\Sl_2(\Z)$ of weight $0$. 
Via the Kronecker limit formula we then study  $I(\tau, \log\|\Delta(z)\|)$. Here  $\|\Delta(z)\|$ is the suitably normalized Petersson metric of the Delta function $\Delta(z)$.

\begin{theorem} \label{thm2:intro}
We have 
\[
I(\tau,\calE_0(z,s)) =  \zeta^{\ast}(s+\frac12) \calF(\tau,s). \tag{i}
\]
Here $\zeta^{\ast}(s)$ is the completed Riemann Zeta function. Moreover,
\[
-\frac1{12} I(\tau, \log\|\Delta(z)\|) =  \calF'(\tau,\frac12), \tag{ii}
\]
where $ \calF'(\tau,\frac12)$ is the derivative of Zagier's Eisenstein series at $s = \tfrac12$.
\end{theorem}

Taking residues at $s=\tfrac12$ in both sides of (i), we obtain another proof that the theta integral~\eqref{thetaintro} is equal to $2\calF(\tau,\frac{1}{2})$. This can be viewed as a special case of the Siegel-Weil formula.

On the other hand, $I(\tau, \log\|\Delta\|)$ can be interpreted in terms of arithmetic geometry. 
In that way, one can recover the main result of \cite{Yang}, to which we refer for 
background information and further details. We let $\mathcal{M}$ be the Deligne-Rapoport compactification of the moduli stack over $\Z$ of elliptic curves, so that 
$\mathcal{M}(\C)$ is the orbifold $\Sl_2(\Z)\back \h \cup \infty$. For $D \in \Z$ and $v>0$, Kudla, Rapoport and Yang \cite{KRY,Yang} construct cycles $\widehat{\calZ}(D,v)$  in the extended arithmetic Chow group of $\calM$ with real coefficients $\cha^1_{\R}(\calM)$,  see \cite{Bost,Kuehn,BKK,Soule}. For $D >0$, the complex points of the underlying divisor of  $\widehat{\calZ}(D,v)$ are the $\G(1)$-equivalence classes of CM points of discriminant $-D$ in $\h$. 
Furthermore, we let $\widehat{\omega}$ be the normalized metrized Hodge bundle on $\calM$, which defines an element $\widehat{c}_1(\widehat{\omega}) = \tfrac{1}{12}( \infty, -\log\|\Delta(z)\|^2)$ in  $\cha_{\R}^1(\calM)$. Finally, we let $\langle \,,\,\rangle$ be the Gillet-Soul\'e intersection pairing. 
Since the divisor of $\Delta$ over $\Z$ does not intersect ${\widehat\calZ}(D,v)$ at the finite places, the $D$-th Fourier coefficient of $-\frac1{12} I(\tau, \log\|\Delta(z)\|)$ turns out to be equal to $4 \langle \widehat{\mathcal{Z}}(D,v), \widehat{\omega}\rangle$.

\begin{theorem} [\cite{Yang}]\label{thm3:intro}
We have 
\begin{equation}\label{ATheta}
\sum_{D \in \Z} \langle \widehat{\mathcal{Z}}(D,v), \widehat{\omega}\rangle q^D =  \frac14 \calF'(\tau,\frac12).
\end{equation}
\end{theorem}

Note that the proof of Theorem~\ref{thm3:intro} of given in \cite{Yang} relies on the explicit calculation and comparison of the Fourier coefficients on both sides of \eqref{ATheta}, while our method does not require that. Also note that we realize the `arithmetic' theta series (Kudla) on the left hand side of \eqref{ATheta} as an honest theta integral. Theorem~\ref{thm3:intro} can be viewed as an instance of an `arithmetic' Siegel-Weil formula envisioned and pursued by Kudla and his collaborators, see e.g.~\cite{Kmsri}, realizing the arithmetic theta series as the derivative of an Eisenstein series.

\medskip

Finally, we show that for $f$ a Maass cusp form of weight $0$, the lift $I(\tau,f)$ is equivalent to a theta lift first introduced by Maass \cite{Maass} and later reconsidered by Katok and Sarnak \cite{KS}.

\medskip

We thank Ulf K\"uhn for suggesting to consider $I(\tau, \log\|\Delta\|)$.
We also thank Gautam Chinta, J\"urg Kramer, Steve Kudla and Steve Rallis for helpful discussions on this project.

\section{Preliminaries}\label{prelim}

Let $V$ be a rational vector space of
dimension $ 3$ with a non-degenerate symmetric bilinear form  $(\,,\,)$ of signature $(1,2)$. We write $q(x) = \tfrac12(x,x)$ for the associated quadratic form and let $d$ be the discriminant of $V$, chosen to be a square-free positive integer.
We fix an orientation for $V$ once and for all. We let $G = \Spin(V) \simeq \SL_2$ viewed as
an algebraic group over $\Q$ and write $\bar{G}\simeq \PSL_2$ for the image in $\Orth(V)$.
We let $D = G(\R)/K $ be the
associated symmetric space, where $K \simeq \SO(2)$ is a maximal compact
subgroup of $G(\R)$. We have $D \simeq \h$, where $\h=\{z\in \C;\;\Im(z)>0\}$ is the complex upper half plane.
For our purposes, it is most convenient to identify $D$ with the space of lines in  $V(\R)$ on which the bilinear form  $(\,,\,)$ is positive definite:
\[
D \simeq \{ z \subset V(\R) ; \;\text{$\dim z =1$ and $(\,,\,)|_z > 0$} \}.
\]

Let $L \subset \V$ be an even lattice of full rank and write $L^\#$ for the dual lattice of $L$. Let $\G$ be a congruence subgroup of $\Spin(L)$ which takes $L$ to itself and acts trivially on the discriminant group $L^{\#}/L$.
We write $M = \G \backslash D$ for the attached locally symmetric space. Throughout we will assume that $M$ is a modular curve, i.e., non-compact. Note that this happens if and only if $V$ is isotropic over $\Q$. We can then view  $V({\Q})$ as the trace zero part $B_0(\Q)$ of the
indefinite quaternion algebra $B(\Q) = M_2(\Q)$. So
\begin{equation} \label{iso}
V(\Q) \simeq  \left\{ X = \begin{pmatrix} x_1 & x_2 \\ x_3 & -x_1
 \end{pmatrix} \in M_2(\Q) \right\}
\end{equation}
with $q(X) = d \det(X)$ and $(X,Y) = -d  \tr(XY)$. In this setting the action of $G \simeq \SL_2$ on $B_0$ is the conjugation:
\[
g.X := gXg^{-1}
\]
for $X\in B_0$ and $g\in {G}$. Moreover, $G(\Q) \simeq \SL_2(\Q)$.

\begin{notation}
From now on, we will write $z= x+iy$ for an
element in the orthogonal symmetric space $D \simeq \h $. The upper case letter $X$  we reserve for vectors in $V(\R)$, thought of as elements in $B_0(\R)$. Its coefficients we denote by $x_i$. Later on, we will write $\tau = u+iv \in \h$ for a modular form variable in $\h$; i.e, we consider $\tau$ as a variable for the (symplectic) symmetric space associated to $\SL_2 \simeq\Symp(1)$.
\end{notation}

We make the previous discussion explicit by giving the
following identification of $D$ with the upper half plane. We pick as base point of ${D}$ the line $z_0$ spanned by $\left( \begin{smallmatrix} 0
& 1 \\ -1 & 0 \end{smallmatrix}\right)$, and note that ${K} = \SO(2)$
is its stabilizer in ${G}(\R)$. For $z \in \h$, we define $g_z \in {G}(\R)/K $ by the condition  $g_zi = z$; the action is the
usual linear fractional transformation on $\h$. We obtain the isomorphism $\h  \to {D}$,
\begin{equation}
z \longmapsto g_z z_0 = \Span\left(g_z . \left( \begin{smallmatrix} 0 & 1 \\ -1 & 0 \end{smallmatrix}\right)\right).
\end{equation}
So for $z = x + iy \in \h$, the associated positive line is generated by
\begin{equation}
X(z) := \frac1{\sqrt{d}} g_z.\zxz{0}{1}{-1}{0} = \frac1{\sqrt{d}y}
\begin{pmatrix} -\frac12(z+ \bar{z}) & z\bar{z} \\ -1 & \frac12(z+ \bar{z})
\end{pmatrix}.
\end{equation}
In particular, $q(X(z)) = 1$ and $g.X(z) = X(gz)$ for $g \in G(\R)$.
For $X= \left( \begin{smallmatrix} x_1 & x_2 \\ x_3 & -x_1
\end{smallmatrix} \right) \in V(\R)$ we have
\begin{align}\label{X,X(z)}
(X,X(z)) &= -\frac{\sqrt{d}}{y} (x_3z\bar{z} - x_1(z+\bar{z}) -x_2) \\
         &=  -\frac{d(x_3x-x_1)^2 +q(X)}{\sqrt{d}x_3y} - \sqrt{d}x_3y,
         \notag
\end{align}
if $x_3\ne0$. We let $(\;,\;)_z$ be the minimal majorant of $(\;,\;)$
associated to $z \in D$. One easily sees that $(X,X)_z = (X,X(z))^2 - (X,X)$.

The set  $\Iso(V)$ of all isotropic lines in $V(\Q)$ can be identified with $P^1(\Q)=\Q\cup\infty$, the set of cusps of $G(\Q)$, by means of the map
\begin{align}\label{def:psi}
\psi:P^1(\Q)\longrightarrow \Iso(V),\quad \psi((\alpha:\beta))=\Span\kzxz{-\alpha\beta}{\alpha^2}{-\beta^2}{\alpha\beta}\in \Iso(V).
\end{align}
One easily checks that $\psi$ is a bijection, commuting with the $G(\Q)$-actions, that is, $\psi(g(\alpha:\beta))=g.  \psi((\alpha:\beta))$. So the cusps of $M$, i.e., the $\Gamma$-classes of  $P^1(\Q)$, can be identified with the $\Gamma$-classes of $\Iso(V)$.
The cusp $\infty\in P^1(\Q)$ is mapped to the isotropic line $\ell_0\in \Iso(V)$ spanned by $X_0= \left( \begin{smallmatrix}  0&1 \\0&0  \end{smallmatrix}\right)$.
For $\ell \in \Iso(V)$, we pick $\sigma_{\ell} \in \SL_2(\Z)$ such that $\sigma_{\ell} \ell_0 = \ell$. We orient all lines $\ell \in \Iso(V)$  by requiring that  $\sigma_{\ell} X_0$ is a positively oriented basis vector of $\ell$.
We let $\G_{\ell}$ be the stabilizer of the line $\ell$. Then (if $-1 \in \G$)
\[
\sigma_{\ell}^{-1}\G_{\ell} \sigma_{\ell} =
 \left\{ \pm
\begin{pmatrix} 1&k\alpha_{\ell} \\0&1 \end{pmatrix} ; \; k \in \Z
 \right\},
\]
where $\alpha_{\ell} \in\Q_{>0}$ is the width of the cusp
${\ell}$. Since $\sigma_{\ell} \in \SL_2(\Z)$, we see that
$\alpha_{\ell}$ does not depend on the choice of $\sigma_{\ell} \in
\SL_2(\Z)$.  For each $\ell$, there is a $\beta_{\ell}\in\Q_{>0}$ such
that $\left( \begin{smallmatrix} 0&\beta_{\ell} \\0 & 0
\end{smallmatrix} \right)$ is a primitive element of $\ell_0 \cap
\sigma_{\ell}^{-1} L$. Finally, we write $\eps_{\ell} =
\alpha_{\ell}/\beta_{\ell}$. Note (see \cite{FCompo}, Definition~3.2)
that $\eps_{\ell}$ would be even well defined if we picked
$\sigma_{\ell} \in \SL_2(\Q)$.  The quantities $\alpha_{\ell}$,
$\beta_{\ell}$, and $\eps_{\ell}$ only depend on the
$\Gamma$-class of $\ell$.

We compactify $M$ to a compact Riemann surface $\bar{M}$ in the usual way by adding a point for each cusp $\ell \in \G \back \Iso(V)$; we also denote this point by $\ell$. For each $\ell \in \Iso(V)$, there is a neighborhood $U_{\ell}$ of $\ell$ such that $z = \g z'$ for some $\g \in \G$ and $z,z' \in U_{\ell}$ implies $\g \in \G_{\ell}$. We write $Q_{\ell} =  e\left (\sigma^{-1}_{\ell} z/\alpha_{\ell}\right)$ with $z \in U_{\ell}$ for the local variable (and for the chart) around $\ell \in \bar{M}$. For $T>0$, we let
$ D_{1/T} = \{w \in \C; \; |w| < \frac{1}{2\pi T}\}$, and note that for $T$ sufficiently big, the inverse images $Q_{\ell}^{-1} D_{1/T}$ are disjoint in $M$.
We truncate $M$ by setting
\begin{equation}\label{truncated}
M_T = \bar{M} - \coprod_{\ell \back \Iso(V)} Q_{\ell}^{-1} D_{1/T}.
\end{equation}

\medskip

In this setting, Heegner points in $M$ are given as follows.
For $X \in V(\Q)$ of positive length, i.e., $q(X) >0$, we put
\begin{equation}
D_X = \Span(X) \in D.
\end{equation}
The stabilizer $G_X$ of $X$ in $G(\R)$ is isomorphic to $\SO(2)$ and for $X \in L^{\#}$, $\G_X = G_X \cap \G$ is finite. We then denote by $Z(X)$ the image of $D_X$ in $M$, counted with multiplicity $\tfrac1{|\bar{\G}_X|}$. We set $D_X = \emptyset$ if $q(X) \leq 0$.

For $m \in \Q_{>0}$ and $h \in L^{\#}$,  $\G$ acts on  ${L}_{h,m} = \{X \in L +h ;\; q(X) =m\}$
 with finitely many orbits. We define the \emph{Heegner divisor} of discriminant $m$  on $M$ by
\begin{equation}
Z(h,m) = \sum_{X \in \G \back L_{h,m} }Z(X).
\end{equation}

On the other hand, a vector $X \in V(\Q)$ of negative length defines a geodesic $c_X$ in $D$ via
\[
c_X = \{ z \in D; \; z \perp X \}.
\]
We denote the quotient $\G_X \back c_X$ in $M$ by $c(X)$. The stabilizer $\bar{\G}_X$ is either trivial (if the orthogonal complement $X^{\perp} \subset V$ is isotropic over $\Q$) or infinite cyclic (if
 $X^{\perp}$ is non-split over $\Q$). If $\G_X$ is infinite, then $c(X)$ is a closed geodesic in $M$, while $c(X)$ is an infinite geodesic if $\bar{\G}_X$ is trivial.
Note that the case $X^{\perp} \subset V(\Q)$ split  is equivalent to $q(X) \in -d \left(\Q^{\times}\right)^2$, see for example  \cite{FCompo}, Lemma~3.6. In that case $X$ is orthogonal to two isotropic lines $\ell_X=\Span(Y)$ and $\tilde{\ell}_X=\Span(\tilde{Y})$, with $Y$ and $\tilde{Y}$ positively oriented. We say $\ell_X$ is the line to associated to $X$ if the triple $(X,Y,\tilde{Y})$ is a positively oriented basis for $V$, and we write $X \sim \ell_X$. Note $\tilde{\ell}_X= \ell_{-X}$.

\section{A Schwartz function of weight $3/2$}

\subsection{Geometric Aspects}

In \cite{KMI}, Kudla and Millson explicitly construct a Schwartz function $\varphi_{KM} = \varphi$ on $V(\R)$ valued in $\Omega^{1,1}(D)$, the differential forms on $D$ of Hodge type $(1,1)$. It is given by
\begin{equation}
\varphi(X,z) = \biggl( (X,X(z))^2  - \frac1{2\pi} \biggr) \, e^{-\pi (X,X)_z}
\, \omega,
\end{equation}
where $ \omega =   \tfrac{dx \wedge dy}{y^2} = \tfrac{i}2 \tfrac{dz \wedge d\bar{z}}{y^2}$. We have $\varphi (g.X,g z) = \varphi (X,z)$ for $g \in G(\R)$.
We define
\begin{align}
\varphi^0(X,z) = e^{\pi(X,X)} \varphi (X,z) =\left((X,X(z))^2  - \frac1{2\pi} \right) \, e^{-2 \pi R(X,z) } \, \omega,
\end{align}
where, following \cite{KAnn}, we set
\begin{equation}\label{Rformel}
R(X,z) := \tfrac12(X,X)_z - \tfrac12(X,X) = \tfrac1{2}(X,X(z))^2 - (X,X).
\end{equation}
The quantity $R(X,z)$ is always non-negative. It equals $0$ if and only if $z =D_X$, i.e., if $X$ lies in the line generated by $X(z)$. Hence, for $X \ne 0$, this does not occur if $q(X) \leq 0$.

The geometric significance of this Schwartz function lies in the fact that for $q(X) >0$, the $2$-form $\varphi^0(X,z)$ is a Poincar\'e dual form for the Heegner point $D_X$,  while $\varphi^0(X,z)$ is exact for $ q(X) <0$. Furthermore,
Kudla \cite{KAnn} constructed a Green function $\xi^0$ associated to  $\varphi^0$. We recall the construction of $\xi^0$.
We consider the exponential integral $\Ei(w)$ for $w \in \C$, defined by $\Ei(w) = \int_{-\infty}^w \tfrac{e^t}{t} dt$, where the path of integration lies in the plane cut along the positive real axis, see e.g. \cite{AbSt}. 
It is well known that $\Ei(w)$ has a logarithmic singularity at $w=0$.
For $X \in V(\R)$, $X \neq0$, we define
\begin{equation}
\xi^0(X,z) = -\Ei(-2\pi R(X,z)).
\end{equation}
Hence  $\xi^0(X,z)$ is a smooth function on $D \setminus
D_X$. For $q(X) >0$, the function $\xi^0(X,z)$ has logarithmic growth at the point $D_X$, while it is smooth on $D$ if $q(X) \leq 0$. In particular,  $\xi^0(X,z)$ is locally integrable. 

We let $\partial$, $\bar{\partial}$ and $d$ be the usual differentials
on $D$. We set $d^c = \tfrac1{4\pi i} (\partial - \bar{\partial})$, so that $dd^c = - \tfrac1{2\pi i} \partial \bar{\partial}$.

\begin{theorem}[Kudla \cite{KAnn}, Proposition 11.1]\label{Green}
Let $X\ne0$. Away from the point $D_X$
\begin{equation}
dd^c \xi^0(X,z) = \varphi^0(X,z).
\end{equation}
The function $\xi^0(X,z)$ is a Green current of logarithmic type for
$D_X$ associated to $\varphi^0(X,z)$ (see \cite{Soule}), i.e., as
currents
\begin{equation}\label{current}
dd^c [\xi^0(X,z)] \; + \; \delta_{D_X} \; = \; [ \varphi^0(X,z)],
\end{equation}
where $\delta_{D_X}$ denotes the delta distribution concentrated at
$D_X$.
\end{theorem}

\begin{proposition}\label{phixigrowth}
For $q(X) >0$, the differential forms $\xi^0(X,z)$, $\partial \xi^0(X,z)$, $\bar\partial \xi^0(X,z)$, and $\varphi^0(X,z)$ are of ``square-exponential'' decay in all directions of $D$, i.e., they are
\begin{align*}
O(e^{-Cx^2}), &\quad \text{as $x \to \pm \infty$,} \\
O(e^{-Cy^2}),& \quad\text{as $y \to \infty$,} \\
O(e^{-C/y^2}),& \quad\text{as $y \to 0$,}
\end{align*}
for some constants $C>0$, and uniformly in $y$ in the first case, and uniformly in $x$ in the other two.
In particular, the current equation (\ref{current}) does not only hold for compactly support functions on $D$, but also for functions of ``linear-exponential'' growth in all directions, i.e., $O(e^{C|x|})$, $O(e^{Cy})$, and $O(e^{C/y})$, respectively.
\end{proposition}

\begin{proof}
Write $X= \left( \begin{smallmatrix} x_1&x_2\\x_3&-x_1 \end{smallmatrix}
\right)$. Since $q(X)>0$, we have $x_3\neq0$. By (\ref{X,X(z)}) we see
\[
{-2 \pi R(X,z)} = { 2 \pi(X,X)} \, - \,{\pi \left(\frac{d(x_3x-x_1)^2 +q(X)}
{\sqrt{d}x_3y} +\sqrt{d}x_3y\right)^2}.
\]
This implies the described decay of the above differential forms.
\end{proof}


\subsection{Automorphic Aspects}




For $\tau = u+ iv \in \h$ we put
$g'_{\tau} =
\left( \begin{smallmatrix}1&u\\0&1\end{smallmatrix} \right)
\left( \begin{smallmatrix}v^{1/2}&0\\0&v^{-1/2}\end{smallmatrix} \right)$,
and define
\begin{align}
\varphi(X,\tau,z) := v^{-3/4} \omega( g'_{\tau}) {\varphi}(X,z)
=  \biggl( v(X,X(z))^2  - \frac1{2\pi} \biggr) \, e^{\pi i (X,X)_{\tau,z}} \,
\omega,
\end{align}
where $(X,X)_{\tau,z} = u(X,X) + iv(X,X)_z = \bar{\tau}(X,X) + iv(X,X(z))^2$. Hence
\begin{equation}\label{phiformel1}
\varphi(X,\tau,z) = e^{2\pi i q(X) \tau} \varphi^0(\sqrt{v}X,z).
\end{equation}

Then, see \cite{KM90, FCompo}, for $h \in L^{\#}/L$, the theta kernel
\begin{equation}
\theta_{h}(\tau,z,\varphi) = \sum_{X \in
h + L} \varphi(X,\tau,z) \in \Omega^{1,1}(D)^{\G}
\end{equation}
 defines a \emph{non-holomorphic}
modular form of weight $3/2$ with values in $\Omega^{1,1}(M)$,
for the congruence subgroup $\G(N)$ of $\SL_2(\Z)$, where $N$ is the level of the lattice $L$ (and for $\G_0(N)$ if $h=0$). More precisely, we let $\Mp_2(\R)$ be the two-fold cover of
$\SL_2(\R)$ realized by the two choices of holomorphic square roots of $\tau \mapsto
j(g,\tau) = c\tau + d$, where $g = \left(
\begin{smallmatrix} a&b \\ c&d \end{smallmatrix} \right) \in
\SL_2(\R)$. Then there is a certain representation $\rho_L$ of the
inverse image $\G'$ of $\Sl_2(\Z)$ in $\Mp_2(\R)$, acting on the group
algebra $\C[L^{\#}/L]$ (see \cite{Bo1}, \cite{Br2}).  We denote the
standard basis elements of $\C[L^{\#}/L]$ by $\frake_{h}$, where $h
\in {L^{\#}/L}$.  For the generators $S = \left( \left(
\begin{smallmatrix} 0&-1\\1&0
\end{smallmatrix} \right), \sqrt{\tau} \right)$,  and $T=\left( \left( \begin{smallmatrix} 1 &1 \\ 0&1 \end{smallmatrix} \right),1\right)$ of $\G'$, the action of $\rho_L$ is given by
\begin{align*}
\rho_L(T) \frake_{h} &= e((h,h)/2) \frake_{h},\\
\rho_L(S) \frake_{h} &= \frac{\sqrt{i}}{\sqrt{|L^{\#}/L|}} \sum_{h'\in L^{\#}/L} e(- (h,h')) \frake_{h'}.
\end{align*}

We then define a vector valued theta series by
\[
\Theta(\tau,z,\varphi) =
 \sum_{h\in L^{\#}/L} \theta_{h}(\tau,z,\varphi) \frake_{h}.
\]
We have, see \cite{KMI,BF},
\[
\Theta(\tau,z,\varphi) \in A_{3/2,L} \otimes \Omega^{1,1}(M),
\]
where $A_{3/2,L}$ denotes the space of $C^{\infty}$-automorphic forms of weight $3/2$ with respect to the representation $\rho_L$, that is, for $(\g',\phi) \in \G'$,
\[
\Theta(\g'\tau,z,\varphi) = \phi^3(\tau) \rho_L(\g',\phi) \Theta(\tau,z,\varphi).
\]
More generally, we denote the holomorphic modular forms of weight $k$ for $\G'$ with respect to $\rho_L$ by $M_{k,L}$, and write $M^!_{k,L}$ for those forms which are holomorphic on $\h$ but meromorphic at the cusp, see e.g. \cite{Bo1,Br2}.

To lighten the notation, we will frequently drop the argument $\varphi$.

\section{The Theta Integral}
\label{sect:4}

We now consider $\Theta(\tau,z)$ as a vector valued top-degree differential form on $M = \G \back D$.
We want to pair it with suitable $0$-forms $f$ on $M$. We need the following result on the growth of $\Theta(\tau,z)$ in $D$.

\begin{proposition}[\cite{FCompo}, Proposition 4.1]\label{thetagrowth}
For each $h \in L^{\#}$ and $\tau \in \h$ and at each cusp $\ell$, we have
\[
\theta_{h}(\tau,\sigma_{\ell} z) = O(e^{-Cy^2}) \qquad \text{as} \qquad y \to \infty,
\]
uniformly in $x$, for some constant $C>0$.
\end{proposition}

\begin{proof}
This follows from the proof of Proposition 4.1, \cite{FCompo}. Note however the confusing typesetting errors in this proof; several occurrences of $\exp(\cdot)$ should be $e(\cdot)$. We therefore give a very brief sketch of the argument given there.

It is easy to see that it is sufficient to assume  $ L = \Z^3 $ in \eqref{iso} and that it suffices to show that $\theta_{h}(\tau,z)$ is rapidly decreasing as $y \to \infty$. For simplicity we assume $d=1$. Note $h= \left( \begin{smallmatrix} h_1 & 0 \\ 0& -h_1 \end{smallmatrix} \right)$ with $h_1 = 0$ or $h_1 = 1/2$. . So we have to consider the growth of
\[
\theta_h(\tau,z)= \sum_{\substack{x_1 \in\Z + h_1 \\ x_2,x_3 \in \Z}} \varphi\left(\left( \begin{smallmatrix} x_1 & x_2 \\ x_3 & -x_1 \end{smallmatrix} \right),\tau,z \right)
\]
as $y \to \infty$. Applying partial Poisson summation with respect to $x_2$, we obtain
\begin{align}\label{thetatrick}
\theta_h(\tau,z) &= -
  \frac{y}{v^{3/2}} \sum_{\substack{x_1 \in\Z + h_1 \\ w,x_3 \in \Z}}
\left(w+x_3\bar{\tau}\right)^2 e\left( -\bar{\tau}x_1^2 \right)
e\left(
-[w+x_3\bar{\tau}][x_3z\bar{z}-2x_1x]\right)  \\ & \quad \times
\exp{\left( -\pi \frac{y^2}{v}(w+x_3\bar{\tau})^2 \right)} dxdy.
\notag \\ &=-
\frac{y}{v^{3/2}} \sum_{\substack{x_1 \in\Z + h_1 \\ w,x_3 \in \Z}}
\left(w+x_3\bar{\tau}\right)^2
e\left( -\bar{\tau}(x_1-x_3x)^2 \right) e\left ( 2(x_1-x_3x/2)xw\right) \notag \\ & \quad \times
\exp\left( -\pi \frac{y^2}{v} |w+ x_3 \tau|^2  \right) dxdy. \notag 
\end{align}
The assertion follows.
\end{proof}

We denote by $ M^!_0(\G)$  the space of (scalar valued) weakly holomorphic modular forms of weight $0$ with respect to $\G$. It consists of those modular functions for  $\G$ which are holomorphic on $D \simeq \h$ and meromorphic at the cusps of $\G$. Hence any $f \in M^!_0(\G)$  has a Fourier expansion at the cusp $\ell$  of the form
\begin{align}\label{fourierf}
f(\sigma_{\ell} z) = \sum_{n \in \tfrac{1}{\alpha_{\ell}} \Z}^{\infty} a_{\ell}(n) e( n z),
\end{align}
with $a_{\ell}(n) =0$ for $n \ll 0$. In particular,
\[
f( \sigma_{\ell} z) = O(e^{2 \pi Ny}) \qquad \qquad (y \to \infty)
\]
for some $N>0$.

We define the theta lift of $f$ by
\begin{align}\label{theta-integral}
I(\tau,f) = \int_{M} f(z) \Theta(\tau,z)
=  \sum_{h \in L^{\#}/L} \left( \int_M f(z) \theta_{h}(\tau,z)  \right) \frake_h.
\end{align}
We also write
\begin{equation}\label{comp-int}
I_h(\tau,f) = \int_M f(z) \theta_{h}(\tau,z)
\end{equation}
for the individual components. Proposition~\ref{thetagrowth} implies the convergence of (\ref{theta-integral}). Then it is clear that $I(\tau,f)$ defines a (in
general non-holomorphic) modular form on the upper half plane of
weight $3/2$.


\begin{definition}(Modular trace for positive index) \\
For $m\in \Q_{>0}$  and $ h \in  L^{\#}/L$, we then define the \emph{modular trace function} of $f$ by
\begin{equation}
{\bf{t}}_f(h,m) = \sum_{z \in Z(h,m)} f(z) = \sum_{X \in \G \back L_{h,m}}
\frac1{|\bar{\G}_X|}
f(D_X).
\end{equation}
\end{definition}

\begin{definition}(Modular trace for $m=0$) \\
For $m=0$, we set
\[
{\bf{t}}_f(h,0) = -  \frac{\delta_{h,0}}{2\pi} \int^{reg}_M f(z) \frac{dx\,dy}{y^2}.
\]
For $f$ non-constant, the integral $\int_M f(z) \frac{dx\,dy}{y^2}$ is divergent, and is regularized by setting
\begin{equation}\label{regularization}
\int^{reg}_M f(z) \frac{dx\,dy}{y^2} = \lim_{T \to \infty} \int_{M_T} f(z) \frac{dx\,dy}{y^2},
\end{equation}
where $M_T$ is the truncated surface defined by \eqref{truncated}.
The regularized integral is computed in Remark \ref{rem:const} below.
\end{definition}

\begin{definition}(Modular trace for negative index) \\
If $n\in\Q_{<0}$ is not of the form $n=-dm^2$ with $m\in \Q_{>0}$ we put
${\bf{t}}_f(h,n)=0$. If $n=-dm^2$ with $m\in \Q_{>0}$ we define
${\bf{t}}_f(h,-dm^2)$ as follows:
Let
$X \in L_{h,-dm^2}$, so that $X^{\perp}$ is split over $\Q$, and $c(X)$ is an infinite geodesic.
We can choose the orientation of $V$ such that
\[
\sigma_{\ell_X}^{-1} X = \left( \begin{matrix} m & r \\ 0 & - m \end{matrix} \right).
\]
 for some $r \in \Q$. In this case the geodesic $c_X$ is explicitly given in $D \simeq \h$ by
\[
c_{X} = \sigma_{\ell_X} \{ z \in \h; \; \Re(z) = -r/2m\}.
\]
We call the quantity $-r/2m$ the {\em real part} of the infinite geodesic $c(X)$ and denote it by $\re(c(X))$. Recall that for the cusp $\ell_X$, we denote the corresponding local variable by $Q_{\ell_X}= e\left(\sigma_{\ell_X}^{-1} z/\alpha_{\ell_X}\right)$. We write $Q_{c(X)} = Q_{\ell_X} e^{2\pi i Re(c(X))/\alpha_{\ell_X}}$.
We now define
\begin{align*}
<f,c(X)> &= - \sum_{n   <0} a_{\ell_X}(n) e^{2\pi i \re(c(X))n } -  \sum_{n <0} a_{\ell_{-X}}(n) e^{2\pi i \re(c(-X))n }\\
&= \Res_{Q_{\ell_X}=0} \left( \frac{f(Q_{c(X)})}{Q_{\ell_X}-1}\right) + \Res_{Q_{\ell_{-X}}=0} \left( \frac{f(Q_{c(-X)})}{Q_{\ell_{-X}}-1}\right) .
\end{align*}
We then put
\[
{\bf{t}}_f(h,-dm^2) = \sum_{X \in \G \back L_{h,-dm^2}} <f,c(X)>.
\]

\end{definition}

\begin{theorem}\label{MAIN}
Let $f \in M_0^!(\G)$ with Fourier expansion as in \eqref{fourierf}, and assume that the constant coefficients of $f$ at all cusps of $M$ vanish. Then the Fourier expansion of $I_h(\tau,f)$ is given by
\begin{align*}
I_h(\tau,f) &=   \sum_{m \geq 0}
{\bf{t}}_f(h,m) q^m  +
 \sum_{\substack{m>0}}    {\bf{t}}_f(h,-dm^2)  q^{-dm^2},
\end{align*}
with $q = e^{2\pi i \tau}$, and where ${\bf{t}}_f(h,m)$ is the modular trace function defined above.

If the constant coefficients of $f$ do not vanish, then $I_h(\tau,f)$ is non-holomorphic, and in the Fourier expansion the following terms occur in addition:
\[
 \frac{1}{2\pi \sqrt{vd}} \sum_{\substack{\ell \in \G \back \Iso(V)\\ \ell \cap L +h \ne \emptyset}} a_{\ell}(0)\eps_{\ell} + \sum_{m>0}
\sum_{X \in \G \back L_{h,-dm^2}}  \frac{a_{\ell_X}(0)+   a_{\ell_{-X}}(0)}{8 \pi \sqrt{vd}m}  \beta(4\pi vdm^2) q^{-dm^2},
\]
where $\beta(s) = \int_1^{\infty}t^{-3/2} e^{-st} dt$.

\end{theorem}

\begin{remark}
(i) The theta lift $I(\tau,f)$ was studied in \cite{FCompo} for the constant function $f =1 \in M_0^!(\G)$. There it was shown that $I_h(\tau,1)$ is {non-holomorphic} and
\[
I_{h}(\tau,1) =
\sum_{m \geq 0} {\bf{t}}_1(h,m) q^m + \frac{1}{2\pi \sqrt{vd}}
\epsilon(h) + \sum_{m>0} \sum_{X \in L_{h,-dm^2}} \frac{1}{4 \pi \sqrt{vd}m}  \beta(4\pi vdm^2) q^{-dm^2}.
\]
Here $\epsilon(h)= \sum_{\ell \in \G \back \Iso(V)} \delta_{\ell}(h) \eps_{\ell}$ with $\delta_{\ell}(h) =1$ if $\ell \cap L +h \ne \emptyset$ and zero otherwise. This generalizes Zagier's non-holomorphic Eisenstein series of weight $3/2$ \cite{Zagier}.

(ii) If $M$ is compact, i.e., a Shimura curve, then $M_0^!(\G) =
M_0(\G) = \C$, and $ I_{h}(\tau,1)$ was considered by Kudla-Millson,
see e.g. \cite{KM90}. Here one has
\[
I_{h}(\tau,1) =
\sum_{m \geq 0} {\bf{t}}_1(h,m) q^m.
\]
\end{remark}

We will now show that the trace function  ${\bf{t}}_f(h,-dm^2)$ vanishes for large $m>0$, so that $I(\tau,f)\in M_{3/2,L}^!$. For this, we sort the infinite geodesics according to the cusps from where they originate.
For $ m \in \Q_{>0}$, we define $L_{h,-dm^2,\ell} = \{ X \in L_{h,-d^2m}; \; X \sim \ell \}$ and see
\[
L_{h,-dm^2} = \coprod_{\ell \in \G \back \Iso(V)} \coprod_{\g \in \G_{\ell} \back \G} \g^{-1} L_{h,-dm^2,\ell}.
\]
Furthermore
\[
\# \G \back L_{h,-dm^2} = \sum_ {\ell \in \G \back \Iso(V)} \# \G_{\ell} \back L_{h,-dm^2,\ell}
\]
so that we conclude
\begin{align}\label{anz}
\nu_{\ell}(h,-dm^2):=\#\G_{\ell} \back L_{h,-dm^2,\ell}  =
\begin{cases}
2 m\eps_{\ell} & \text{if $L_{h,-dm^2,\ell} \ne \emptyset$} \\
0 & \text{else}.
\end{cases}
\end{align}
with $\eps_{\ell} = \alpha_{\ell}/\beta_{\ell}$ as in section~\ref{prelim} (see \cite{FCompo}, Lemma~3.7).

\begin{proposition}\label{prop:neg}
Let $f \in M_0^!(\G)$ with Fourier expansion as in \eqref{fourierf}. Then
\begin{align*}
{\bf{t}}_f(h,-dm^2) &= -   \sum_{\ell \in \G \back \Iso(V)}  \nu_{\ell}(h,-dm^2) \sum_{n \in \frac{2m}{\beta_{\ell}} \Z_{ <0}} a_{\ell}(n) e^{2\pi i  r n }   \\
&\quad -   \sum_{\ell \in \G \back \Iso(V)}  \nu_{\ell}(-h,-dm^2)
\sum_{n \in \frac{2m}{\beta_{\ell}}\Z_{<0} }
a_{\ell}(n) e^{2\pi i  r' n },
\end{align*}
with $r=\re(c(X))$ for any $X \in L_{h,-dm^2,\ell}$ and $r'=\re(c(X))$ for any $X \in L_{-h,-dm^2,\ell}$.
In particular,
\[
{\bf{t}}_f(h,-dm^2) = 0 \qquad \qquad \text{for} \qquad  m \gg 0.
\]

\end{proposition}


\begin{proof}
We have
\begin{align}\label{oben}
{\bf{t}}_f(h,-dm^2) &= - \sum_{X \in \G \back L_{h,-dm^2}}   \sum_{n \in \frac{1}{\alpha_{\ell}} \Z_{<0}}      a_{\ell_X}(n) e^{2\pi i \re(c(X))n } \\ &\quad -   \sum_{X \in \G \back L_{h,-dm^2}}   \sum_{n \in \frac{1}{\alpha_{\ell}}  \Z_{<0}} a_{\ell_{-X}}(n) e^{2\pi i \re(c(-X))n }. \notag
\end{align}
We denote the first term in \eqref{oben} by $G(h,-dm^2)$.
So the second term in \eqref{oben} is equal to $G(-h,-dm^2)$. We have
\[
G(h,-dm^2) = \sum_{\ell \in \Gamma \back \Iso(V)} G(h,-dm^2,\ell),
\]
where
\[
G(h,-dm^2,\ell) =  - \sum_{X \in \G \back L_{h,-dm^2,\ell}}  \sum_{n \in \frac{1}{\alpha_{\ell}} \Z_{<0}} a_{\ell}(n) e^{2\pi i \re(c(X))n }
\]
We can assume that a set of representatives for $\G_{\ell} \back L_{h,-dm^2,\ell}$ is given by
\begin{equation*}\label{reps}
\left\{  Y_k = \sigma_{\ell} \,m\zxz{1}{2r + k\beta_{\ell}/m}{0}{-1}; \quad k = 0,\dots, 2m \eps_{\ell} -1 \right\}
\end{equation*}
for some $r\in \Q$. In particular, $\re(c(Y_k)) = -r - k\frac{\beta_{\ell}}{2m}$. Thus,
\begin{align*}
G(h,-dm^2,\ell) &= -\sum_{k=0}^{2\eps_\ell m -1}\sum_{n  \in \Z_{<0}} a_{\ell}(n/\alpha_{\ell}) e^{-2\pi i (r + k\beta_{\ell}/2m  )n/\alpha_{\ell} } \\
 &=-\sum_{n  \in \Z_{<0}}  a_{\ell}(n/\alpha_{\ell})  e^{-2\pi i rn/\alpha_{\ell}}  \sum_{k=0}^{2\eps_\ell m -1} e^{-2\pi i n  k/(2m\eps_{\ell})} \\
&= - 2m\eps_{\ell}\sum_{n  \in 2m\eps_{\ell} \Z_{<0}}  a_{\ell}(n/\alpha_{\ell})  e^{-2\pi i rn/\alpha_{\ell}}.
\end{align*}
The other term,  $G(-h,-dm^2)$, is treated in the same way.
\end{proof}

Theorem~\ref{MAIN} and Proposition~\ref{prop:neg} imply

\begin{corollary}
Assume that all constant coefficients of $f \in M_0^!(\G)$ vanish. Then
\[
I(\tau,f)\in M_{3/2,L}^!.
\]
\end{corollary}

\begin{remark}\label{rem:const}
One can compute ${\bf{t}}_f(h,0)$ as follows. We consider the Eisenstein series for $\G$ of weight $2$ at the cusp
$\ell_0$, i.e., at $\infty$:
\begin{equation}
E_2(z,s)= \sum_{\g \in \G_{\ell_0} \back \G} j(\g,z)^{-2}|j(\g,z)|^{-2s},
\end{equation}
where $ j(g,z)= {cz+d}$ for $g = \left(
\begin{smallmatrix} a &b \\ c&d
\end{smallmatrix} \right)$.
Then, see e.g., \cite{Kubota}, 
the series $E_2(z,s)$ converges for $s >0$ and has a
meromorphic continuation to $\C$. At $s=0$, $E_2(z,s)$ is
holomorphic, and we put $E_2(z) = E_2(z,0)$ which defines a (non-holomorphic) modular form of weight $2$ for $\G$. The Fourier expansion $E_{2,\ell}(z)= j(\sigma_{\ell},z)^{-2} E_2(\sigma_{\ell} z)$ at a cusp $\ell$ is of the form
\begin{equation}
E_{2,\ell}(z) = \left(b_\ell(0) +
c(0)\frac1{y}\right) + \sum_{n=1}^{\infty} b_{\ell}(n/ \alpha_{\ell}  )
e^{2\pi i nz/\alpha_{\ell}}.
\end{equation}
Here $b_\ell(0)=\delta_{\ell,\ell_0}$ is the Kronecker delta and $c(0)=-\frac{\alpha_{\ell_0}}{\vol_\omega(M)}$ is independent of $\ell$.

Using Stokes' theorem and the fact that $\bar\partial (E_{2,\ell}(z)dz)=c(0)\frac{dx\,dy}{y^2}$  one sees similarly to \cite{Bo1} section 9 that the regularized divergent integral $\int^{reg}_{\G \back D} f(z) d\mu$ is equal to
\begin{align}
-\frac{1}{c(0)}\sum_{\ell \in \G \back \Iso(V)}   \alpha_\ell
\sum_{\substack{n\in \frac{1}{\alpha_\ell}\Z_{\geq 0}}} a_\ell(-n)b_\ell(n).
\end{align}

In particular, if $\G$ is a congruence subgroup of $\SL_2(\Z)$, we may make this more explicit, using the Fourier expansion
\[
\calE_2(z)=-\frac{3}{\pi y}-24\sum_{n=0}^\infty \sigma_1(n)e^{2\pi i nz}
\]
of the (non-holomorphic) Eisenstein series $\calE_2(z)$ of weight $2$ for $\Sl_2(\Z)$. Here $\sigma_1(0)=-\frac{1}{24}$ and $\sigma_1(n)=\sum_{t\mid n} t$ for $n\in \Z_{>0}$.
Arguing as in \cite{Bo1} section 9 we find in this case
\begin{align}
\int^{reg}_{\G \back D} f(z) \frac{dx\,dy}{y^2}
=-8\pi \sum_{\ell \in \G \back \Iso(V)}   \alpha_\ell
\sum_{\substack{n\in \Z_{\geq 0}}} a_\ell(-n)\sigma_1(n).
\end{align}
For $f=1$ we recover the well known relation $\sum_{\ell \in \G \back \Iso(V)}   \alpha_\ell =[\PSL_2(\Z):\overline{\Gamma}]$.
\end{remark}

\begin{proof}[Proof of Theorem~\ref{MAIN}]
We give the outline of the structure of the proof, which reduces the theorem to the computation of several orbital integrals. We will compute these integrals in the next section.

We define
\begin{equation}\label{Fourier1}
\theta_{h,m}(\tau,z) = \sum_{X \in {L}_{h,m}} \varphi(X,\tau,z)  \qquad \text{and}
\qquad \theta^0_{h,m}(v,z) = \sum_{X \in {L}_{h,m}}
\varphi^0(\sqrt{v}X,z).
\end{equation}
By \eqref{phiformel1} we then have
\begin{align}\label{Fourierexp}
I_{h}(\tau,f) = \int_{M} \sum_{m\in\Q}
f (z) \theta_{h,m}(\tau,z)  = \sum_{m\in\Q} \left( \int_{M} f (z)
\theta^0_{h,m}(v,z) \right) q^m,
\end{align}
which is the Fourier expansion of $I_{h}(\tau,f)$. (Hence interchanging summation and integration is valid in the last step).

For $m \ne 0$, $\G \back {L}_{h,m}$ is finite. Therefore, for these $m$, we obtain for the latter integral in (\ref{Fourierexp}):
\begin{align}\label{Fourier+}
\int_{M} f (z) \theta^0_{h,m}(v,z) &=
\int_{M} \sum_{X\in \G \backslash {L}_{h,m}  } \, \sum_{\g \in
\G_X \backslash \G} f (z) \varphi^0(\g^{-1}\sqrt{v}X,z) \\
&= \sum_{X\in \G \backslash    {L}_{h,m} } \, \int_{M} \sum_{\g \in \G_X \backslash \G} f (z)
\varphi^0(\sqrt{v}X,\g z), \notag
\end{align}
provided the interchange of summation and integration is valid, i.e., the integral
in (\ref{Fourier+}) converges for all $X$.

Then the statement about the positive Fourier coefficients of $I(\tau,f)$ follows from

\begin{proposition}\label{prop1}
Let $X \in L +h$ such that $q(X) >0$.
Then
\[
\sum_{\g \in \G_X \backslash \G} f (z)
\varphi^0(\sqrt{v}X,\g z) \in L^1(M),
\]
and
\[
 \int_{M} \sum_{\g \in \G_X \backslash \G} f (z)
\varphi^0(\sqrt{v}X,\g z) = \frac1{|\bar{\G}_X|} f(D_X).
\]
\end{proposition}


\medskip

For $q(X) < 0$, the space $X^{\perp}\subset V$ has signature $(1,1)$, and we have to  distinguish two cases, depending on whether  $X^{\perp}$ is isotropic over $\Q$ or not. If $X^{\perp}$ is isotropic over $\Q$, then $\bar{\G}_X$ is trivial and $q(X) \in -d \left(\Q^{\times}\right)^2$. If not, $\bar{\G}_X$ is infinite cyclic  and $q(X) \notin  -d \left(\Q^{\times}\right)^2$ (see \cite{FCompo} Lemma 4.2).

For $m \notin - d \left(\Q^{\times}\right)^2$, \eqref{Fourier+} reduces the statement about the $m$-th coefficient to

\begin{proposition}\label{prop2a}
Let $X \in L_{h,m}$ with $m <0$ such that $m \notin - d \left(\Q^{\times}\right)^2$. Hence $\bar{\G}_X$ is infinite cyclic.
Then
\[
\sum_{\g \in \G_X \backslash \G} f(z)
\varphi^0(\sqrt{v}X, \g z) \in L^1(M)
\]
and
\begin{equation*}
 \int_{M} \sum_{\g \in \G_X \backslash \G} f(z)
\varphi^0(\sqrt{v}X, \g z) = 0.
\end{equation*}
\end{proposition}


For the split case, we have

\begin{proposition}\label{prop2b}
Let  $X \in L_{h,-dm^2}$ (with $m \in \Q_{>0}$) so that $\bar{\G}_X=1$. Then
\[
 \int_{M} \sum_{ \g \in \G } f(z)
\varphi^0(\sqrt{v}X, \g z) \in L^1(M)
\]
and
\begin{align*}\label{prop2bformel}
 \int_{M} \sum_{ \g \in \G}  f(z) \varphi^0(\sqrt{v}X, \g z)&= \left(a_{\ell_X}(0)+   a_{\ell_{-X}}(0) \right)\frac{1}{8 \pi \sqrt{vd}m}  \beta(4\pi vdm^2)
 \\ &\quad -  \sum_{n  <0} a_{\ell_X}(n) e^{2\pi i \re(c(X))n }  -     \sum_{n <0} a_{\ell_{-X}}(n) e^{2\pi i \re(c(-X))n } .
 \end{align*}
\end{proposition}


\medskip

It remains to compute the constant coefficient of  $I_h(\tau,f)$, which is given by
\begin{equation}\label{isoint}
\int_M \sum_{\substack{X \in L + h \\ q(X) =0} } f(z) \varphi^0(\sqrt{v}X,z).
\end{equation}
We would like to split this integral into two pieces; one for $X=0$ (if $h=0$) and the other for $X \ne 0$. However, for $X=0$, we simply have $\varphi(\tau,X) = -\tfrac1{2\pi}\omega$ and therefore $\int_{M} f(z) \varphi^0(0,z) $ does \emph{not} converge due to the exponential growth of $f$. In order to split the integral (\ref{isoint}) we therefore have to regularize it, as explained in \eqref{regularization}. We obtain
\begin{equation}\label{regulsplit}
\int_M \theta_{h,0}^0(\sqrt{v},z) = -\frac{\delta_{h,0}}{2\pi} \int_M^{reg} f(z) \, \omega +
\int_M^{reg} \sum_{\substack{X \in L_{h,0} \\ X \ne 0 } } f (z) \varphi^0(\sqrt{v}X,z).
\end{equation}
The first term is ${\bf{t}}_f(h,0)$.

\begin{proposition}\label{prop3}
For the second regularized integral in \eqref{regulsplit}, we have
\[
\int_M^{reg} \sum_{\substack{X \in L_{h,0}\\ X \ne 0 } } f (z)\varphi^0(\sqrt{v}X,z)  = \frac{1}{2\pi\sqrt{vd}} \sum_{\substack{\ell \in \G \back \Iso(V)\\ \ell \cap L+h \ne \emptyset}}  a_{\ell}(0) \eps_{\ell}.
\]
\end{proposition}



This finishes the (outline of the) proof of Theorem~\ref{MAIN}.
\end{proof}

\section{Orbital Integrals}

In this part of the paper, we will prove Propositions~\ref{prop1}, \ref{prop2a}, \ref{prop2b}, and \ref{prop3}.

We begin with a lemma on Fourier transforms, which we will need later. For a function $g(t)$ on the real line, let $\hat{g}(w)= \int_{-\infty}^{\infty} g(t) e^{2\pi itw} dt$ be its Fourier transform.

\begin{lemma}\label{Ftransform}
For $a,b >0$, let
\[
h(t) =  t(t -ib) \frac{e^{-a^2 t^2}}{t^2 + b^2}.
\]
Then
\[
\hat{h}(w) = -\frac{1}{a} e^{a^2 b^2} \left(\pi ab e^{2\pi bw} \erfc\left(\tfrac{a^2b+\pi w}{a}\right) - \sqrt{\pi} e^{ - \tfrac{a^4b^2 +\pi^2  w^2}{a^2}} \right).
\]
Here $\erfc(x)$ is the standard complementary error function given by
\[
\erfc(x) = \frac{2}{\sqrt{\pi}} \int_x^{\infty} e^{-u^2}du.
\]
\end{lemma}

\begin{proof}[Proof of Lemma~\ref{Ftransform}]
By \cite{E}, p.~74 (26), the Fourier transform of $ f(t) = t \frac{e^{-a^2 t^2}}{t^2 + b^2}$
is
\[
\hat f(w) = \frac{ \pi i }2 e^{a^2b^2} \left( e^{-2\pi bw} \erfc(ab - \pi w/a) - e^{2\pi bw} \erfc(ab + \pi w/a) \right),
\]
(note the different normalization there). By differentiating under the integral, we see that the Fourier transform of $ tf(t) = t^2 \frac{e^{-a^2 t^2}}{t^2 + b^2}$ is given by the derivative $-\frac{i}{2 \pi} \hat f'(w)$. Since $h(t) = tf(t) - ibf(t)$, we obtain for the Fourier transform of $h$:
\[
\hat h(w)=-\frac{i}{2 \pi} \hat f'(w)-ib \hat f(w).
\]
But
\begin{align*}
-\frac{i}{2 \pi} \hat f'(w) & =  -\frac{\pi b}{2} e^{a^2b^2}\left(
 e^{-2\pi bw} \erfc(ab - \pi w/a)  + e^{2\pi bw} \erfc(ab + \pi w/a)
\right) \\ & \quad+ \frac{2 \sqrt{\pi}}{a}  e^{a^2b^2} e^{-(a^2b^2 + \pi^2w^2/a^2)}.
\end{align*}
Lemma~\ref{Ftransform} follows.
\end{proof}

\begin{proof}[Proof of Proposition~\ref{prop1}]
Let $X\in L^{\#}$ such that $q(X)>0$. Then $\bar{\G}_X$ is a finite cyclic group. Using the $\G$-invariance of $f$, we see
\begin{align}\label{prop1cal}
\int_{\G \back D} f (z)\sum_{\g \in\G_X \backslash \G}  \varphi^0(\sqrt{v} X, \g z) & = \int_{ \G_X \back D} f (z) \varphi^0(\sqrt{v}X, z) \\
\notag
& = \frac1{|\bar{\G}_X|} \int_D  f(z)\varphi^0(\sqrt{v}X, z).
\end{align}
By Proposition~\ref{phixigrowth}, the decay of $\varphi^0(\sqrt{v}X,z)$ offsets the growth of $f$. Therefore the last integral in (\ref{prop1cal}) exists, which implies the existence of the first integral
and the validity of the unfolding. By Theorem~\ref{Green}, Proposition~\ref{phixigrowth}, and $D_{\sqrt{v}X} = D_X$ we see
\[
\int_D f(z) \varphi^0(\sqrt{v}X,z) =  f(D_X) + \int_D  \xi^0(\sqrt{v}X,z)dd^c f (z).
\]
But $dd^c f  =0$, since $f$ is holomorphic. This proves Proposition~\ref{prop1}.
\end{proof}

\begin{proof}[Proof of Proposition~\ref{prop2a}]

Let $q(X)=m<0$ for $X\in V$, so $X^{\perp}$ has signature $(1,1)$. Assume that $X^{\perp}$ is non-split, so that $\G_X$ is infinite cyclic.
By conjugation, we can always assume that $X = \sqrt{-m/d}
\left( \begin{smallmatrix}
1&0\\0&-1 \end{smallmatrix} \right)$. Then
$\G_X=  \left\langle \left( \begin{smallmatrix}
\epsilon&0\\0&\epsilon^{-1} \end{smallmatrix} \right) \right\rangle $
with some $\epsilon >1$.
Using \eqref{X,X(z)} we find for our particular choice of $X$ that  $(X,X(z))^2 = -4m\tfrac{x^2}{y^2}$.
Therefore, in view of the explicit formula for $\varphi^0(X,z)$, we obtain by (formally) unfolding the orbital integral:
\begin{align*}
\int_{\G \backslash D} f(z) \sum_{\g \in \G_X \backslash \G}
\varphi^0(\sqrt{v}X, \g z) &= \int_{\G_X \backslash D} f(z) \varphi^0(\sqrt{v}X, z)  \\
&= e^{4\pi mv} \int_{\G_X \backslash D}
f (z) \left( -4mv \frac{x^2}{y^2} - \frac{1}{2\pi} \right) e^{4\pi mv \frac{x^2}{y^2}} \frac{dx \, dy}{y^2}.
\end{align*}
A fundamental domain $\mathcal{G}$ of $\G_X \backslash D$ is the domain bounded by the semi arcs $|z|=1$ and $|z|= \epsilon^2>1$ in the upper half plane:
\begin{equation}
\mathcal{G} = \left\{ z \in D ; \; 1 \leq |z| < \epsilon^2 \right\}.
\end{equation}
But in this region, the rapid decay of $\varphi^0(X,z)$ offsets the
growth of $f(z)$ as $z$ approaches the boundary of $\mathcal{G}$. So all
considered integrals actually exist and unfolding is allowed. Finally,
by Theorem~\ref{Green} we have
\[
\int_{\G_X \backslash D} f(z) \varphi^0(\sqrt{v}X, z) = \int_{\G_X \backslash D}  \xi^0(\sqrt{v}X, z) dd^c f(z) =0,
\]
since $f$ is holomorphic. This proves Proposition~\ref{prop2a}.
\end{proof}

\begin{proof}[Proof of Proposition~\ref{prop2b}]
Here we consider the case that $ q(X) = -d m^2$ ($m>0$). Note that the proof of Proposition~\ref{prop2a} does not carry over, since for $X\in L_{h,-dm^2}$  and $\bar{\G}_X$ trivial, the integral $ \int_{D} f(z) \varphi^0(\sqrt{v}X, z) $ does \emph{not} exist. (Even for $f=1$, see \cite{FCompo}).
Since $f$ is holomorphic, by Stokes' theorem we have
\begin{align*}
\int_{M} f(z) \sum_{\g \in \G} \varphi^0(\sqrt{v}X, \g z)
&=
\frac1{2\pi i} \int_{M} f (z) \bar{\partial} \partial \sum_{\g \in \G} \xi^0(\sqrt{v}X, \g z) \\
&=
\frac1{2\pi i} \int_{M } d \left( f(z)  \partial \sum_{\g \in  \G} \xi^0(\sqrt{v}X, \g z) \right) \\
& = \frac1{2\pi i} \lim_{T \to \infty} \int_{\partial M_T }
 f (z)  \sum_{\g \in \G}  \partial \xi^0(\sqrt{v} X, \g z).
\end{align*}
Note here
\begin{equation}\label{dxiformel}
\partial \xi^0(X, z) = - \frac{\partial R(X,  z)}{R(X,z)} e^{-2\pi R(X,z)}.
\end{equation}
For an isotropic line $\ell$, we write $\partial M_{T,\ell}$ for the boundary component of $M_T$ at the cusp corresponding to $\ell$. So $\partial M_T = \coprod_{\ell \in \G \back \Iso(V)} \partial M_{T,\ell}$.

For any $X \in L_{-dm^2}$, there is an involution $J_X \in G(\Q)$ taking $X$ to $-X$ and interchanging the lines $\ell_X$ and $\tilde{\ell}_X$. (It could be made unique be requiring in addition that $J_X (\ell_X \cap L) = \tilde{\ell}_X \cap L$.)
For example, for $X = m \left( \begin{smallmatrix} 1 & 2r \\ 0 &-1 \end{smallmatrix} \right)$, we can take $J_X = T_{-r} J T_{r}$ where $J= \left(\begin{smallmatrix} 0&1 \\ -1&0 \end{smallmatrix} \right)$ and $T_r = \left(\begin{smallmatrix} 1&r \\ 0&1 \end{smallmatrix} \right)$.
So for an arbitrary $X \in L_{-dm^2}$, we can take $J_X = \sigma_{\ell_X} J_{X'} \sigma_{\ell_X}^{-1}$ where $X' = \sigma_{\ell_X}^{-1}X$.


\begin{lemma}\label{lemmagisela}

\begin{align}\label{gisela}
\lim_{T \to \infty} \int_{\partial M_T  }
 f (z)  \sum_{\g \in \G}  \partial \xi^0(\sqrt{v} X, \g z)
 &= \lim_{T \to \infty}  \int_{\partial M_{T,\ell_X}  }
 f ( z) \sum_{ \g \in \G_{\ell_X}} \partial \xi^0(\sqrt{v} X,  \g  z) \\
\nonumber
&\phantom{=}{}+
\lim_{T \to \infty}   \int_{\partial M_{T,\tilde{\ell}_X}  }
 f (z)  \sum_{ \g \in \G_{\tilde{\ell}_X}} \partial \xi^0(\sqrt{v} X,    \g z).
\end{align}

\end{lemma}

\begin{proof}
Choosing the orientation of $V$ appropriately, we have
\[
X':= \sigma_{\ell_X}^{-1} X = m \begin{pmatrix} 1 & 2r \\ 0 &-1 \end{pmatrix}
\]
for some $r \in \Q$. Then
\begin{align}\label{gisbert}
\lim_{T \to \infty} & \int_{\partial M_T} f(z) \sum_{\g \in \G}  \partial \xi^0(\sqrt{v}X, \g z) \notag \\
&= -\lim_{T \to \infty} \sum_{\ell \in \G \back \Iso(V)} \int_{z=iT}^{\alpha_{\ell}+iT}
 f (\sigma_{\ell} z) \sum_{\g \in \G} \partial \xi^0(\sqrt{v} X, \g\sigma_{\ell} z) \\
 &= - \lim_{T \to \infty} \sum_{\ell \in \G \back \Iso(V)} \int_{z=iT}^{\alpha_{\ell}+iT}
 f (\sigma_{\ell} z) \sum_{\g \in \G} \partial \xi^0(\sqrt{v} X', \sigma_{\ell_X}^{-1}\g\sigma_{\ell} z). \notag
\end{align}

We have
\begin{align}\label{Xformula}
(X',X(z))^2 = 4d m^2 \frac{(x+ r ) ^2}{y^2} &=4d m^2\left(\frac{1}{\Im(z)\Im(J(z+r))}-1\right).
\end{align}
If $g =\kabcd\in G(\R)$, we see by means of \eqref{Rformel} and \eqref{Xformula} that
\begin{align*}
R(X', g z) = 2dm^2 \frac{1}{\Im(gz)\Im({J_{X'} gz})}
=2dm^2 \frac{|cz+d|^2|(a+rc)z+b+rd|^2}{y^2}
\end{align*}
with $J_{X'} = T_{-r}JT_r$ as above.

Let $\tilde\Gamma$ be an arithmetic subgroup of $G(\Q)$. Then
there is an $\eps>0$ such that  $R(X', g z)>\eps$ for all
$g\in \tilde \Gamma$, uniformly on $y>1$.
Moreover, using \eqref{dxiformel}, one easily checks that there is a $\delta>0$ such that
\begin{align*}
\partial \xi^0(\sqrt{v}X',g z)
&\ll e^{-\delta (|cz+d|^2+|(a+rc)z+b+rd|^2)}e^{-\delta y^2} dz
\end{align*}
for all $g=\kabcd\in \tilde\Gamma$
with $c\neq 0$ and $a+rc \neq 0$, uniformly for $y>1$.


This implies that
\[
-  \lim_{T \to \infty} \sum_{\ell \in \G \back \Iso(V)} \int_{z=iT}^{\alpha_{\ell}+iT}
 f (\sigma_{\ell} z) \sum_{\substack{\g \in  \G\\c( \sigma_{\ell_X}^{-1}\g \sigma_{\ell})\neq 0, \\ c( J_{X'} \sigma_{\ell_X}^{-1}\g \sigma_{\ell}) \ne 0}}   \partial \xi^0(\sqrt{v} X',  \sigma_{\ell_X}^{-1} \g \sigma_{\ell} z)=0,
\]
where $c(g)$  denotes the lower left entry of $g \in \Sl_2(\R)$.
Consequently, in (\ref{gisbert}) we only have to consider the terms with   $ c( \sigma_{\ell_X}^{-1}\g \sigma_{\ell}) =0$ or $ c( J_{X'} \sigma_{\ell_X}^{-1}\g \sigma_{\ell})  =0$.
But  $ c( \sigma_{\ell_X}^{-1}\g \sigma_{\ell}) =0$ is equivalent to $   \sigma_{\ell_X}^{-1}\g \sigma_{\ell}  \ell_0 = \ell_0$. Hence $\ell =\sigma_{\ell} \ell_0$ is $\G$-equivalent to $\ell_X =\sigma_{\ell_X} \ell_0$ and therefore we may assume $\ell= \ell_X$.  Now   $ c( \sigma_{\ell_X}^{-1}\g \sigma_{\ell_X}) =0$  implies $\g \in \G_{\ell_X}$.
We obtain the first summand on the right hand side of (\ref{gisela}).

On the other hand,   $ c( J_{X'} \sigma_{\ell_X}^{-1}\g \sigma_{\ell})  =0$   means $\g \sigma_{\ell} \ell_0 = \sigma_{\ell_X} J_{X'} \ell_0= J_X {\ell}_X  = \tilde{\ell}_X$. Hence $\ell$ is $\G$-equivalent to $\tilde{\ell}_X$. So we assume $\ell =  \tilde{\ell}_X$ and hence $\g \in  \G_{\tilde{\ell}_X}$. This gives rise to the second summand on the right hand side of (\ref{gisela}).
\end{proof}

\begin{lemma}\label{britta}
For $X \in L_{h,-dm^2}$, we have
\begin{align}
&\frac1{2\pi i} \lim_{T \to \infty}  \int_{\partial M_{T,\ell_X}  }
 f ( z) \sum_{ \g \in \G_{\ell_X}} \partial \xi^0(\sqrt{v} X,  \g  z) \\
\nonumber
&=\frac{1}{8\pi \sqrt{vd}m  } a_{\ell_X}(0) \beta(4\pi vdm^2) - \sum_{n \in \frac1{\alpha_{\ell_X}}\Z_{<0}} a_{\ell_X}( n)  e^{2 \pi i  \re(c(X))   n}.
\end{align}

\end{lemma}

\begin{proof}
As before, we can write $X':= \sigma_{\ell_X}^{-1} X = m\left( \begin{smallmatrix} 1 & 2r \\ 0 &-1 \end{smallmatrix} \right)$ for some $r \in \Q$. Hence $\re(c(X)) = -r$. For simplicity, we write $\alpha= \alpha_{\ell_X}$ and $g(z)= f(\sigma_{\ell_X} z)$ with Fourier expansion $g(z)=\sum_{n \in \frac{1}{\alpha} \Z} a(n) e(nz)$. We first see
\begin{align}\label{limit}
&\frac1{2\pi i} \lim_{T \to \infty}  \int_{\partial M_{T,\ell_{X}}  }
 f ( z) \sum_{ \g \in \G_{\ell_X}} \partial \xi^0(\sqrt{v} X,  \g  z) \\ &= -\frac1{2\pi i}
\lim_{T \to \infty}  \int\limits_{z=iT}^{\alpha+iT}
 g ( z) \sum_{n \in \Z} \partial \xi^0\left(\sqrt{v} m
\left( \begin{smallmatrix}  1& 2(r+\alpha n) \\0 & -1  \end{smallmatrix} \right),  z\right)  \nonumber
\end{align}
%
For $Y= m  \left( \begin{smallmatrix}  1&  2(r+\alpha n)  \\0 & -1 \end{smallmatrix} \right)$, we note
\begin{align*}
R(Y,z) &= \frac{2dm^2}{y^2}(x+ r + \alpha n)^2 +2dm^2,          \\ 
\partial R(Y,  z)
&= \frac{2dm^2}{y^2}(x+ r + \alpha n  ) (1 +\frac{i}{y}(x +r +  \alpha n   )).
\end{align*}
Therefore by \eqref{dxiformel}:
\begin{multline*}
\partial \xi^0( \sqrt{v} Y, z)   =  \frac{-i}{y}e^{ -4\pi dm^2 v} (x + r+ \alpha n)  ( x + r+ \alpha n -i  y)
\frac{e^{ -4 \pi v dm^2 (x + r+ \alpha n)^2/y^2}}{(x + r+ \alpha n)^2 +  y^2} \, dz.
\end{multline*}
We set $t = x + r+ \alpha n$,
$a = 2 \frac{\sqrt{\pi vd}m}{y}$, and
$b = y$, and obtain
\[
\partial \xi^0( \sqrt{v} Y, z) = -\frac{ i}{b}e^{- a^2 b^2}
t(t -ib) \frac{e^{-a^2 t^2}}{t^2 + b^2} \, dz = -\frac{ i}{b}e^{- a^2 b^2} h(t) dz
\]
with $h(t)$ as in Lemma~\ref{Ftransform}. Hence, the Fourier transform of $h_1(t) = h(x+ r+ \alpha t)$ is given by
\[
\hat{h}_1(w) = \frac{1}{ \alpha} e^{-2\pi i (x+r) w/\alpha} \hat h(w/\alpha).
\]
By Poisson summation,  Lemma~\ref{Ftransform} therefore gives
\begin{equation*}
\begin{split}
&\sum_{n\in \Z} \partial \xi^0\left( \sqrt{v}m \left( \begin{smallmatrix}  1& 2(r + \alpha n) \\0 & -1 \end{smallmatrix} \right), z\right)  \\ &=
\sum_{w \in \frac1{\alpha} \Z} \frac{i}{2\alpha \sqrt{ vd}m  }   e^{-2\pi i (x+r) w}  \\ &  \qquad \times
\left( 2 \pi \sqrt{vd}m e^{2\pi wy } \erfc\left( 2\sqrt{\pi vd}m + \sqrt{\pi} wy/ 2\sqrt{vd}m \right) -  e^{-4\pi vdm^2- \pi  w^2y^2/ 4vdm^2} \right) dz
 \end{split}
\end{equation*}
Inserting the Fourier expansion for $g$ and carrying out the integration
we get for the quantity in \eqref{limit}:
\begin{align*}
&-\frac{1}{4\pi \sqrt{ vd}m  }\lim_{T \to \infty}  \sum_{w \in \frac1{\alpha} \Z}
 a( w) e^{-2\pi i rw} e^{-2\pi wT }\\
& \quad \times\left( 2 \pi \sqrt{vd}m e^{2\pi wT} \,\erfc\left( 2\sqrt{\pi vd}m + \sqrt{\pi} wT/ 2 \sqrt{vd }m\right) -  e^{-4\pi vdm^2 -\pi  w^2T^2/ 4 vdm^2} \right) \\
& =  -\frac{1}{4\pi\sqrt{ vd}m  }\lim_{T \to \infty}  \sum_{w \in \frac1{\alpha} \Z}
 a( w) e^{-2\pi i r w}\\
& \quad \times
\left( 2 \pi \sqrt{vd}m \, \erfc\left( 2\sqrt{\pi vd}m + \sqrt{\pi} wT/ 2\sqrt{vd }m\right) -
e^{-\pi(2\sqrt{ vd}m + wT/ 2\sqrt{vd}m)^2} \right)
.
\end{align*}
The square exponential decay of $e^{-\pi(2\sqrt{ vd}m + wT/ 2\sqrt{vd}m)^2}$ for $w\neq 0$ implies that the contribution corresponding to these terms vanishes in the limit.
Therefore the above quantity is equal to
\begin{align}\label{jochen}
& -\frac{1}{4\pi\sqrt{ vd}m  }
 a(0) \left( 2 \pi \sqrt{vd}m \, \erfc\left( 2\sqrt{\pi vd}m \right) -  e^{-4\pi vdm^2} \right)\\
\nonumber
&-\frac{1 }{2 }\lim_{T \to \infty}  \sum_{w \in \frac1{\alpha}\Z\setminus\{0\}}
 a( w)  e^{-2\pi i rw} \, \erfc\left( 2\sqrt{\pi vd}m + \sqrt{\pi} wT/ 2\sqrt{vd }m\right).
\end{align}
Using the identity $\beta(t)= 2\left(e^{-t} - \sqrt{\pi t} \erfc(\sqrt{t})\right)$ we find that the first term in \eqref{jochen} is equal to
\[
  \frac{1}{8\pi\sqrt{ vd}m}a(0) \beta(4\pi vdm^2).
\]
For the second term in \eqref{jochen}, we first note that $\erfc(t) = O(e^{-t^2})$ as $t\to + \infty$ and $\lim_{t \to - \infty} \erfc{t} =2$. Hence
 the second term in \eqref{jochen} is equal to
\[
-\sum_{w \in \frac1{\alpha}\Z_{<0}}
 a( w)  e^{-2\pi i rw}.
\]
This gives  Lemma \ref{britta}.
\end{proof}

This finishes the proof of Proposition~\ref{prop2b}.
\end{proof}

\begin{proof}[Proof of Proposition~\ref{prop3}]

We now consider Proposition~\ref{prop3}; the sum over the non-zero isotropic vectors. We write $X_{\ell}$ for the primitive positive oriented vector in $L \cap \ell$. We can write $ \ell \cap (L+h) = \Z X_{\ell} +h_{\ell}$ for some $h_{\ell} \in L+h$ if $\delta_{\ell}(h) \ne 0$. We then have
\begin{align*} 
\int_{\G \backslash D}^{reg} f(z) \sum_{\substack{X\in L_{h,0} \\ X\neq0}}
\varphi^0(\sqrt{v}X,z)
&=\int_{\G \backslash D}^{reg} f(z) \sum_{\ell \in \G \back \Iso(V)}
\sum_{\substack{X\in \ell \cap (L+h)\\ X\neq0}} \sum_{ \g \in \G_{\ell} \back \G} \varphi^0(\sqrt{v} \g^{-1} X,z) \\
&=      \sum_{\substack{\ell \in \G \back \Iso(V) \\ \delta_{\ell}(h) \ne 0}}     \int_{\G \backslash D}^{reg} f(z)
 \sum_{\g\in\G_{\ell} \back \G} \sideset{}{'}{\sum}_{n=-\infty}^{\infty}
\varphi^0(\sqrt{v}(nX_{\ell}+h_{\ell}), \g z). \notag
\end{align*}
Here $\sum'$ indicates that we omit $n=0$ in the sum in the case of the trivial coset. As before, we obtain by Stokes' theorem
\begin{multline}\label{formeliso}
\int_{\G \backslash D}^{reg} f(z) \sum_{\substack{X\in L_{h,0} \\ X\neq0}}
\varphi^0(\sqrt{v}X,z) \\ =  \frac1{2\pi i}     \sum_{\substack{\ell \in \G \back \Iso(V) \\ \delta_{\ell}(h) \ne 0}}      \lim_{T \to \infty}
\int_{\partial M_T} f(z)
\sum_{\g\in\G_{\ell} \back \G} \sideset{}{'}{\sum}_{n=-\infty}^{\infty}
\partial \xi^0(\sqrt{v}(nX_{\ell}+h_{\ell}), \g z).
\end{multline}
Note $(X, X(z)) = \sqrt{d}r/y$ for $X = \left(\begin{smallmatrix} 0&r \\ 0&0\end{smallmatrix} \right)$.  By \eqref{dxiformel} we find for
$ g=\kabcd \in G(\R)$ that
\[
\partial \xi^0(\sqrt{v}X, g z) =  -\frac{i}{(cz+d)^2 \Im(gz)}
e^{-\pi v d r^2/\Im(gz)^2} dz.
\]
Similarly to the proof of Proposition~\ref{prop2b}, we then see that on the right hand side of \eqref{formeliso} in the limit the terms for $\g  \ne 1$ vanish, while for $\g =1$, we have a contribution at the boundary component corresponding to the cusp $\ell$. Thus
\begin{equation*}
\int_{\G \backslash D}^{reg} f(z) \sum_{\substack{X\in L_{h,0} \\ X\neq0}}
\varphi^0(\sqrt{v}X,z) =
\frac1{2\pi}
\sum_{\substack{\ell \in \G \back \Iso(V) \\ \delta_{\ell}(h) \ne 0}}     \lim_{T \to \infty}
\int_{iT}^{iT + \alpha_{\ell}}\!\!\!  f(z) \sideset{}{'}{\sum}_{n=-\infty}^{\infty} \frac1y e^{-\pi vd (n\beta_{\ell}+k_{\ell})^2 /y^2} dx.
\end{equation*}
Here $\sigma_{\ell}^{-1} (X_{\ell}+ h_{\ell}) = \left(\begin{smallmatrix} 0&\beta_{\ell} + k_{\ell} \\ 0&0\end{smallmatrix} \right)$ for some number $k_{\ell}$.
Note that in the limit a possible term for $n=0$ and $k_{\ell}=0$ vanishes. Then, by carrying out the integral and Poisson summation we obtain
\begin{align*}
\int_{\G \backslash D}^{reg} f(z) \sum_{\substack{X\in L_{h,0} \\ X\neq0}}
\varphi^0(\sqrt{v}X,z) &=     \sum_{\substack{\ell \in \G \back \Iso(V) \\ \delta_{\ell}(h) \ne 0}}  \frac{\alpha_{\ell}}{2\pi} a_{\ell}(0) \lim_{T \to \infty}\sum_{n=-\infty}^{\infty} \frac1T e^{-\pi vd (n\beta_{\ell}+k_{\ell})^2 /T^2} \\
 & =  \sum_{\substack{\ell \in \G \back \Iso(V) \\ \delta_{\ell}(h) \ne 0}} \frac{\eps_{\ell}a_{\ell}(0) }{2\pi\sqrt{vd}}     \lim_{T \to \infty} \sum_{w \in \Z}   e(-w k_{\ell}/\beta_{\ell}) e^{-\pi w^2 T^2/(vd\beta_{\ell}^2)}\\
 & =   \sum_{\substack{\ell \in \G \back \Iso(V) \\ \delta_{\ell}(h) \ne 0}}   \frac{\eps_{\ell}}{2\pi\sqrt{vd}} a_{\ell}(0).
\end{align*}
This concludes the proof of Proposition~\ref{prop3}.
\end{proof}

\section{Example}\label{sec:ex}

We explain how to obtain the example from the introduction.
Let $p$ be a prime. We consider the quadratic space $V(\Q)$ as in \eqref{iso} with the quadratic form $q(X)=\det(X)$. We let $L$ be the lattice
\[
L=\left\{\zxz{b}{2c}{2ap}{-b};\quad a,b,c\in \Z\right\}.
\]
Then $L$ has level $4p$ and is stabilized by $\G_0(p)$.
The modular curve $M=\Gamma_0(p)\bs D$ is compactified by adding the two cusps $\infty,0$ of $\Gamma_0(p)$, which are represented by the isotropic lines
\begin{equation}\label{isolines}
\ell_0=\Span \zxz{0}{1}{0}{0},\qquad \ell_1=\Span\zxz{0}{0}{-1}{0}.
\end{equation}
We may take $\sigma_{\ell_0}=1$ and $\sigma_{\ell_1}=\kzxz{0}{-1}{1}{0}$. One checks that $\alpha_{\ell_0}=1$, $\beta_{\ell_0}=2$, $\eps_{\ell_0}=1/2$, and
$\alpha_{\ell_1}=p$, $\beta_{\ell_1}=2p$, $\eps_{\ell_1}=1/2$.

The Heegner points now can be described as follows. If $X=\kzxz{b}{2c}{-2ap}{-b}\in L$ is a vector of positive norm $-\Delta=q(X)$, then the matrix
\begin{equation}
Q= \zxz{ap}{b/2}{b/2}{c}= \frac{1}{2}\zxz{0}{-1}{1}{0}X
\end{equation}
defines a definite integral binary quadratic form of discriminant $\Delta=b^2-4pac=-q(X)$. Here the $\Gamma_0(p)$-action on $L$ corresponds to the natural right action on quadratic forms, and the cycle $D_X$ coincides with the CM point $\alpha_{Q}$ (resp. $\alpha_{-Q}$) corresponding to $Q$ (resp. $-Q$) if $Q$ is positive (resp. negative) definite as in the introduction.
We then easily see
\begin{equation}\label{AZ}
Z(0, -\Delta)= \sum_{Q\in \calQ_{-\Delta,p}/\Gamma_0(p)}\frac{2}{|\Gamma_0(p)_Q|}\alpha_Q.
\end{equation}

Let $f\in M_0^!(\Gamma_0(p))$ be a weakly holomorphic modular form and denote its Fourier expansions at the cusps $\infty$, $0$ by
\begin{align*}
f(z)&=\sum_{n\in \Z}a(n)e(nz) \qquad \text{and} \qquad 
f(\sigma_{\ell_1} z)=\sum_{n\in \frac{1}{p}\Z}b(n)e(nz),
\end{align*}
respectively. By \eqref{AZ}, we have
\begin{equation}\label{A+}
{\bf{t}}_f(0, -\Delta)
= \sum_{Q\in \calQ_{-\Delta,p}/\Gamma_0(p)}\frac{2}{|\G_0(p)_Q|}f(\alpha_Q).
\end{equation}
By means of Remark \ref{rem:const}, we see that
\begin{equation*}
{\bf{t}}_f(0, 0)=4 \sum_{\substack{n\in \Z_{\geq 0}}} \big( a(-n)\sigma_1(n)+ p b(-n)\sigma_1(n)\big).
\end{equation*}
We find a different expression for ${\bf{t}}_f(0, 0)$ by applying the residue theorem to the meromorphic $1$-form $f(z)\big(\calE_2(z)-\calE_2|_2(W_p) (z)  \big) dz$ on $\overline{\G_0(p) \back \h}$. This yields
\[
 \sum_{\substack{n\in \Z_{\geq 0}}}  a(-n) \big(\sigma_1(n)- p \sigma_1(n/p)\big) =  \sum_{\substack{n\in \Z_{\geq 0}}}  b(-n/p) \big(\sigma_1(n)- p \sigma_1(n/p)\big),
\]
and therefore
\begin{equation}\label{A0}
{\bf{t}}_f(0, 0)= 2 \sum_{\substack{n\in \Z_{\geq 0}}}  \big( a(-n) + b(-n/p) \big) \big(\sigma_1(n)+ p \sigma_1(n/p)\big).
\end{equation}

For the modular traces of $f$ with negative index $n$, we first recall that by Proposition \ref{prop:neg}, we have ${\bf{t}}_f(0,n)=0$ unless $n=-m^2$ with $m\in \N$. Furthermore, $\kzxz{m}{0}{0}{-m}\in L_{0,-m^2,\ell_0}$ and $\kzxz{-m}{0}{0}{m}\in L_{0,-m^2,\ell_1}$. This implies that the quantities $r$ and $r'$ in Proposition \ref{prop:neg} are equal to $0$. Thus
\begin{equation}\label{A-}
{\bf{t}}_f(0, -m^2)= -2m\sum_{k\in \Z_{>0}} \big( a(-mk)+b(-mk/p) \big).
\end{equation}

Collecting the terms \eqref{A+}, \eqref{A0}, \eqref{A-} now shows that
Theorem \ref{MAIN} implies  Theorem \ref{thm:intro} of the introduction: For $f \in M^!_0(\G^{\ast}_0(p))$ (i.e., $f$ is in the $+1$-eigenspace for the Fricke involution $W_p$), we have $a(n) = b(n/p)$, and 
${\bf{t}}_f(0, N) = 2 {\bf{t}}^{\ast}_f(N)$ for $N>0$. Thus, if $a(0)=0$, then
\[
G(\tau,f)=\frac{1}{4} I_0(\tau,f).
\]
Finally note that $-q(X)$ is congruent to a square modulo $4p$ for $X \in L$ (which we write as $-q(X)\equiv \square\;(4p)$). Consequently, $G(\tau,f)$ belongs to $M^{+,!}_{3/2}(p)$, the Kohnen plus space of weakly holomorphic modular forms of weight $3/2$ for the group $\Gamma_0(4p)$ having a Fourier expansion of the form
\begin{align}\label{kohnencond}
g(\tau)=\sum_{\substack{n\in \Z\\-n\equiv \square\;(4p)}} c(n)q^n.
\end{align}

If $f \in M^!_0(\G_0(p))$ is in the  $-1$-eigenspace for $W_p$, we have 
 $a(n) = -b(n/p)$, and 
\[
 I_0(\tau,f) =0,
\]
since we directly see ${\bf{t}}_f(0, N) =0$ for $N>0$, while for $N\leq 0$ we have
${\bf{t}}_f(0, N) =0$ by  \eqref{A0}, \eqref{A-}.

For $p=1$, we get $G(\tau,f)=\frac{1}{2} I_0(\tau,f)$, and for $f=J$, we recover Zagier's result.

\section{Extensions}

In this section, we consider other automorphic forms of weight $0$ for $\G$ as input for the theta lift under consideration in this paper.

\subsection{The Lift of the weight $0$ Eisenstein Series and $\log|\Delta|$}

For $z \in \h$ and $s \in \C$, we let
\[
\calE_0(z,s) = \frac12 \zeta^{\ast}(2s+1) \sum_{\g \in \G_{\infty} \back \Sl_2(\Z)} \left(\Im(\g z)\right)^{s+\frac12}
\]
be the (normalized) real analytic Eisenstein series of weight $0$ for $\Sl_{2}(\Z)$. Here $\G_{\infty} = \left( \begin{smallmatrix}1&\Z \\0&1 \end{smallmatrix} \right)$ and $\zeta^{\ast}(s) = \pi^{-s/2} \G(\tfrac{s}2) \zeta(s)$ is the completed Riemann Zeta function. Recall that 
$\calE_0(z,s)$ converges for $\Re(s) >1/2$ and has a meromorphic continuation to $\C$ with a simple pole at $s = 1/2$ with residue $1/2$. Furthermore, it is well known that $\calE_0(z,-s) = \calE_0(z,s)$.


We  consider the quadratic space $V(\Q)$ as in \eqref{iso} with the quadratic form $q(X)=\det(X)$. For simplicity, we let $L$ in this section be the lattice
\[
L=\left\{\zxz{b}{c}{a}{-b};\quad a,b,c\in \Z\right\}.
\]
We have $L^\#/L\cong \Z/2\Z$, the level of $L$ is $4$, and $\Gamma=\Sl_2(\Z)$ takes $L$ to itself and acts trivially on $L^\#/L$. We let $\mathfrak{e}_0, \mathfrak{e}_1$
be the standard basis of  $\C[L^\#/L]$ corresponding to the cosets $h= \left( \begin{smallmatrix} h_1 & 0 \\ 0& -h_1 \end{smallmatrix} \right)$ with $h_1 = 0$ and  $h_1 = 1/2$, respectively.

We let $K$ be the one-dimensional lattice $\Z$ together with the negative definite bilinear form $(b,b') = -2 bb'$. We naturally have $L^{\#}/L \simeq K^{\#} / K$.
We define a vector valued Eisenstein series $\calE_{3/2,K}(\tau,s)$ of weight $3/2$ for the representation $\rho_{K}$ by
\[
\calE_{3/2,K}(\tau,s) = -\frac1{4\pi} (s+\frac12) \zeta^{\ast}(2s+1) \sum_{\g' \in \G'_{\infty} \back \G'} \left( v^{\frac12(s-\frac12)} \mathfrak{e}_0\right)|_{3/2,K} \,\g',
\]
where the Petersson slash operator is defined on functions $f: \h \to \C[K^{\#}/K]$ by
\[
(f |_{3/2,K}\, \g') (\tau) =  \phi(\tau)^{-3} \rho^{-1}_{K}(\g') f(\g \tau)
\]
for $\g'=(\g,\phi)\in \G'$.
Here $\G_{\infty}'$ is the inverse image of $\G_{\infty}$ inside $\G'$.
Again we have $\calE_{3/2,K}(\tau,-s)= \calE_{3/2,K}(\tau,s)$, as we will also see below. We set
\begin{equation}\label{ZEisen}
\mathcal{F}(\tau,s) =  \left(\calE_{3/2,K}(4\tau,s) \right)_{0} + \left( \calE_{3/2,K}(4\tau,s) \right)_{1} .
\end{equation}
Then the value of $\mathcal{F}(\tau,s)$ at $s= 1/2$ is a (non-holomorphic) modular form of weight $3/2$  for $\G_0(4)$ and is equal to Zagier's Eisenstein series as in \cite{HZ,Zagier}. This can be seen as follows. The right hand side of \eqref{ZEisen} realizes the isomorphism of vector valued modular forms of type $\rho_K$ with the space of modular forms for $\G_0(4)$ satisfying the Kohnen plus-space condition, see \cite{EZ}, section~5. On the other hand, Zagier's Eisenstein series is the only Eisenstein series of weight $3/2$  for $\G_0(4)$ in the plus-space and has the same constant coefficient as $\mathcal{F}(\tau,\tfrac12)$. Note that our $\mathcal{F}(\tau,s)$ has a different normalization as Zagier's  $\mathcal{F}(\tau,s)$, see also \cite{Yang}, section 3.

\begin{theorem}\label{Eisensteinlift}
With the notation as above, we have
\begin{equation}\label{th71}
I(\tau,\calE_0(z,s)) =  \zeta^{\ast}(s+\frac12) \calE_{3/2,K}(\tau,s).
\end{equation}
\end{theorem}

\begin{proof}

 As in \cite{Bo1}, section~4 we define the theta series
\[
\Theta_{K}(\tau, \alpha,\beta) = \sum_{h \in  K^{\#} / K} 
\sum_{ x_1 \in K + h} e\left( - \bar{\tau} (x_1+\beta)^2\right) e\left( -(x_1 +\beta/2, \alpha)  \right) \mathfrak{e}_h.
\]
By \eqref{thetatrick} we then have
\begin{align*}
\Theta&(\tau,z)  = -\frac{y}{v^{3/2}}\sum_{w, x_3 \in\Z} 
\left(w+x_3\bar{\tau}\right)^2 \exp\left( -\pi \frac{y^2}{v} |w+ x_3 \tau|^2  \right) \Theta_{K}(\tau, -wx,-x_3x) dxdy\\
& =  -\frac{y}{v^{3/2}} \sum_{n=1}^{\infty} n^2\sum_{\substack{c, d \in\Z \\ gcd(c,d) =1} }
\left(c\bar{\tau} +d \right)^2 \exp\left( -\pi \frac{n^2y^2}{v} |c\tau +d|^2 \right)
 \Theta_{K}(\tau, -ndx,-ncx)dxdy.
 \notag 
\end{align*}
Now take $a,b \in \Z$ such that $\g' =  \left( \kabcd , \sqrt{c\tau+d}\right)\in \Gamma'$. By \cite{Bo1}, Theorem~4.1 we find
\begin{equation}\label{Borchformula}
 \Theta_{K}(\tau, -ndx,-ncx) = \left(c\bar{\tau} +d \right)^{-1/2} \rho_{K}^{-1} \left( \g' \right)\Theta_{K}(\g'\tau, -nx,0).
\end{equation}
Hence
\begin{multline*}
\Theta(\tau,z)  =  -\frac{y}{v^{3/2}} \sum_{n=1}^{\infty} n^2\sum_{\g'\in \G'_{\infty} \back \G'}
\left(c\bar{\tau} +d \right)^{3/2} \exp\left( -\pi \frac{n^2y^2}{v} |c\tau +d|^2 \right)
\\ \times \rho_{K}^{-1}\left( \g'\right) \Theta_{K}(\g'\tau, -nx,0)dxdy.
\end{multline*}
Then by the standard Rankin-Selberg unfolding trick we obtain for $\Re(s)>1$:
\begin{align*}
I(\tau,\calE_0(z,s)) & = \zeta^{\ast}(2s+1)\int_{ \G_{\infty} \back \h} \Theta(\tau,z) y^{s+\frac12} \\
& =   -{v^{-3/2}} \zeta^{\ast}(2s+1) \sum_{n=1}^{\infty} n^2 
\sum_{\g' \in \G'_{\infty} \back \G'} \left(c\bar{\tau} +d \right)^{3/2}  \\
& \quad \times  \int_0^{\infty}  \exp\left( -\pi \frac{n^2y^2}{v} |c\tau +d|^2 \right) y^{s+\frac52} \frac{dy}{y} \\
& \quad \times
\rho_{K}^{-1}\left( \g' \right)  \left( \int_0^1  \Theta_{K}(\g'\tau, -nx,0) dx \right) \\
& = -\zeta^{\ast}(2s+1) \G\left(\frac12(s+\frac12)+1 \right) \pi^{-(1+\frac12(s+\frac12))}  \zeta(s+\frac12)  \\ 
& \quad \times \frac12
\sum_{\g' \in \G'_{\infty} \back \G'} \frac{v^{\frac12(s-\frac12)}}{ |c\tau+d|^{s-\frac12}} \frac1{\left(c{\tau} +d \right)^{3/2}} \rho_{K}^{-1}\left( \g' \right) \mathfrak{e}_0 \\
& = \zeta^{\ast}(s+\frac12) \calE_{3/2,K}(\tau,s).
\end{align*}
\end{proof}

Taking residues at $s=1/2$ on both sides of \eqref{th71} we obtain

\begin{corollary}\label{ZagierSW}
\[
I(\tau,1) =  2 \calE_{3/2,K}\left(\tau,\frac12\right).
\]
\end{corollary}


We let 
\[
\Delta(z) = e^{2\pi i z} \prod_{n=1}^{\infty} \left(1-e^{2\pi i nz}\right)^{24}
\]
 be the Delta function. We normalize the Petersson metric of $\Delta$ such that
\[
\|\Delta(z)\| = e^{-6C} |\Delta(z) (4\pi y)^6|,
\]
with $C = \tfrac12 (\g + \log4\pi)$.

\begin{theorem}
We have
\[
- \frac{1}{12} I\left(\tau,\log\left( \|\Delta(z)\| \right)\right) = \calE'_{3/2,K}\left(\tau,\frac12\right).
\]

\end{theorem}

\begin{proof}

Recall that the Kronecker limit formula states
\begin{equation}\label{Kroni}
- \frac{1}{12}\log\left( |\Delta(z)y^6| \right) = \lim_{s \to \frac12} \left( \calE_0(z,s) - \zeta^{\ast}(2s-1)\right).
\end{equation}
By \eqref{Kroni}, Theorem~\ref{Eisensteinlift} and Corollary~\ref{ZagierSW} we have
\begin{align*}
- \frac{1}{12}I(\tau,\log\left( |\Delta(z)y^6| \right)) &=  \lim_{s\to\frac12} 
\Bigl( I\left(\tau, \calE_0(z,s) \right)- I(\tau,\zeta^{\ast}(2s-1)) \Bigr)\\
& =  \lim_{s\to\frac12}   \Bigl(\zeta^{\ast}(s+\frac12) \calE_{3/2,K}(\tau,s) -   2 \zeta^{\ast}(2s-1)   \calE_{3/2,K}(\tau,\frac12) \Bigr) \\
& =  \calE'_{3/2,}(\tau,\frac12) + \frac12 \left(\log(4\pi)-\g \right) \calE_{3/2,K}(\tau,\frac12).
\end{align*}
Here we used
\[
\zeta^{\ast}(s) =  \frac1{s-1} - \frac12(\log (4\pi)-\g) + O(s-1).
\]
The theorem follows now follows from 
\[
-\frac{1}{2}I(\tau,\log(4\pi)-C) = - 
 \frac12 \left(\log(4\pi) -\g \right) \calE_{3/2,K}(\tau,\frac12).
\]
\end{proof}

With the notation as in the introduction, the cycles $ \widehat{\mathcal{Z}}(m,v)$ for $m>0$ with $-m \equiv 0,1 \mod 4$, are given by, see \cite{Yang} section~3,
\[
 \widehat{\mathcal{Z}}(m,v) = (\mathcal{Z}(m), \Xi(m,v)) \in  \cha_{\R}^1(\calM).
\]
Here $\mathcal{Z}(m)$ is the divisor in $\calM$ given by the moduli stack over $\Z$ of elliptic curves $E$ such that there is an embedding $\mathcal{O}_{m} \hookrightarrow \End(E)$, where $\mathcal{O}_m$ is the order of discriminant $-m$ in $\Q(\sqrt{-m})$. Thus $\mathcal{Z}(m)(\C) = \calQ_{m}/\Sl_2(\Z)=: Z(m)$ (with each elliptic curve counted with multiplicity $\tfrac{1}{\#\Aut (E)}$). Moreover,
\[
\Xi(m,v) = \frac14   \sum_{\substack{X \in L^{\#} \\ \frac12(X,X)=m}}  \xi^0(2\sqrt{v}X)
\]
is a Green function for $Z(m)$. For $m \leq 0$, the  $ \widehat{\mathcal{Z}}(m,v)$ are defined similarly using $\xi^0$ with the divisor either supported at $\infty$ (if $m=-n^2$ or $m=0$) or empty (otherwise). 

\begin{remark}
In \cite{Yang} and in section~\ref{sec:ex} (with $p=1$), the cycles 
$Z(m)$ are constructed using the trivial coset of the 
lattice $\tilde{L} = \{\left( \begin{smallmatrix} b & 2c \\ 2a & -b \end{smallmatrix} \right);
\; a,b,c\in \Z \}$ in $V(\Q)$. Since $\tilde{L} = 2L^{\#}$, we can use $L$ instead. On the other hand, for the proof of Theorem~\ref{Eisensteinlift}, the setting of vector valued modular forms and theta series, in particular \eqref{Borchformula}, is quite convenient. Via \eqref{ZEisen} we then can go back to the scalar valued situation.
\end{remark}

\begin{theorem}
We have
\begin{equation}
\sum_{m \in \Z} \langle \widehat{\mathcal{Z}}(m,v), \widehat{\omega}\rangle q^m =  \frac14 \calF'(\tau,\frac12).
\end{equation}

\end{theorem}

\begin{proof}
We only show this for $m>0$. For the other coefficients we refer to \cite{Yang}; they can be done with the methods developed in this paper as well. We have
\[
\calF'(\tau,\frac12) = - \frac{1}{12} \int_M \sum_{X \in L^{\#}}  \varphi^0(2\sqrt{v}X,Z)
\log\left( \|\Delta(z)\| \right) 
 e^{\pi i (X,X) 4 \tau}.
\]
This follows from \eqref{th71}. 
By \eqref{current} we have
\begin{align}\label{arithint}
- \frac{1}{12} \sum_{\substack{X \in L^{\#} \\ \frac12(X,X)=m}}  \int_M \varphi^0 (2\sqrt{v}X) 
\log\left( \|\Delta(z)\| \right) &=  - \frac{1}{3} \sum_{z \in Z(m) } \log\left( \|\Delta(z)\| \right) \\& \quad + \frac{1}{2\pi} \int_M \Xi(m,v) \frac{dxdy}{y^2}.\notag
\end{align}
Since the divisor of $\Delta$ over $\Z$ is disjoint to $\calZ(m)$, we now easily see using the definition of the star product that \eqref{arithint} is equal to  $4 \langle \widehat{\mathcal{Z}}(m,v), \widehat{\omega}\rangle$.
\end{proof}

Our method should generalize to modular curves of higher level. Furthermore, the results above suggest that one should consider $I(\tau,\log\|f\|)$ for other modular forms than $\Delta$. In particular, the case when $f$ is a Borcherds lift \cite{Bo1,Br2} could be of interest.

\subsection{The lift of Maass cusp forms}

We let $L_{cusp}^2(\G \back D)$ be the space of cuspidal square integrable functions on $\G \back D = M$. It is clear that we consider $I(\tau,f)$ for $f \in  L_{cusp}^2(\G \back D)$ as well. It turns out that this lift is closely relating to another theta lift first considered by Maass \cite{Maass} and later reconsidered by Katok and Sarnak \cite{KS}. Namely, they considered, in our notation, the space $V^-$, which is the space $V$ together with the negative bilinear form $-(\,,\,)$. Hence $V^-$ has signature $(2,1)$. The Siegel theta series for $V^-$ is given by
\[
\theta_h(\tau,z,\varphi_{2,1}) = \sum_{X \in L +h} \varphi_{2,1}(X,\tau,z)
\]
with $\varphi_{2,1}(X,\tau,z) = v e^{\pi i(-u(X,X) + iv(X,X)_z)}$. Then  $ \theta_h(\tau,z,\varphi_{2,1})$ is automorphic with weight $1/2$ for $\tau \in \h$. We can then define
\[
I_M(\tau,f) = \sum_{h \in L^{\ast}/L} \left( \int_M f(z) \theta_h(\tau,z,\varphi_{2,1})  \frac{dx\,dy}{y^2} \right) \frake_h
\]
for $f \in L_{cusp}^2(\G \back D)$. In fact, in \cite{Maass,KS} only Maass forms are considered, that is, eigenfunctions of the hyperbolic Laplacian $\Delta = -y^2\left(  \frac{\partial^2}{\partial x^2}  +  \frac{\partial^2}{\partial y^2}\right)$.

For the relationship between $I$ and $I_M$, first recall that the Maass raising and lowering operators are given by $
R_k = 2i \frac{\partial}{\partial \tau} + k v^{-1}$ and $L_k = -2iv^2  \frac{\partial}{\partial \bar{\tau}}$.
Hence $R_{k-2}L_k = -\Delta'_k$, where  $\Delta'_k$
is the weight $k$ Laplacian for $\tau \in \h$ as in \cite{Br2}.
We also need the operator $\xi_k$ which maps forms of weight $k$ to forms of  ``dual'' weight $2-k$. It is given by
\[
\xi_k(f)(\tau) = v^{k-2} \overline{L_k f(\tau)} = R_{-k} v^k \overline{f(\tau)}.
\]
\begin{lemma}
The two kernel functions $\varphi_{KM} =\varphi$ and $ \varphi_{2,1}$ of the two lifts $I$ and $I_M$ satisfy the following fundamental relationship:
\[
\xi_{1/2} \varphi_{2,1}(X,\tau,z) \cdot\omega= - \pi \varphi_{KM}(X,\tau,z).
\]
Furthermore, we have
\[
-4 \Delta'_{\frac12} \varphi_{2,1}(X,\tau,z) =  \Delta \varphi_{2,1}(X,\tau,z).
\]
\end{lemma}
\begin{proof}
This can be easily seen by a direct and straightforward calculation. Alternatively, one can switch to the Fock model of the Weil representation, see e.g. \cite{BF}, section 4, and perform the calculation there.
\end{proof}


\begin{theorem}
For $f \in  L_{cusp}^2(\G \back D)$, we have
\[
\xi_{1/2} I_M(\tau,f) =  - \pi I(\tau,f).
\]
If $f$ is an eigenfunction of $\Delta$ with eigenvalue $\la$, then we also have
\[
\xi_{3/2}  I(\tau,f) =  -\frac{\la} {4\pi} I_M(\tau,f).
\]
\end{theorem}

\begin{proof}
The first assertion immediately follows from the lemma. For the second, note that
$\xi_{3/2} \xi_{1/2} = R_{-3/2} L_{1/2} = - \Delta'_{1/2}$. Then by the adjointness of $\Delta$ we see
\begin{align*}
\xi_{3/2}  I(\tau,f)   = - \frac{1}{\pi} \xi_{3/2} \xi_{1/2} I_M(\tau,f) =
\frac{1}{\pi} \Delta_{1/2}' I_M(\tau,f)  
 =   -\frac1{4\pi} I_M(\tau,\Delta f) =
-\frac{ \la}{4\pi} I_M(\tau,f).
\end{align*}
\end{proof}

The theorem shows that the two lifts are equivalent on Maass forms. Note however, that due to the moderate growth of $\theta_h(\tau,z,\varphi_{2,1})$ one cannot define $I_M(f)$ on $M_0^{!}$. On the other hand, since $I(\tau,f)$ is holomorphic for $f \in M_0^!$, we have $\xi_{3/2}  I(\tau,f) =0 $.

\subsection{The lift of weak Maass forms}

In \cite{BF}, section 3, we introduced the space of weak Maass forms $H_k(\G)$. It consists of those forms $f(z)$ on $D$ of weight $k$ for $\G$ which are annihilated by the weight $k$ Laplacian and satisfy $f(\sigma_{\ell}z ) = O(e^{Cy})$ as $z \to \infty$ for some constant $C$. Here we are only interested in $H_0(\G)$. A form $f \in H_0(\G)$ can be written as $f = f^+ + f^-$, where the Fourier expansions of $f^+$ and $f^-$  are of the form
\begin{align*}
f^+(\sigma_{\ell} z) &= \sum_{n \in \frac1{\alpha_{\ell}} \Z} a^+_{\ell} (n) e(nz) \\
f^-(\sigma_{\ell} z) &= a_{\ell}^-(0) v + \sum_{\substack{n \in \frac1{\alpha_{\ell}} \Z-\{0\}}} a^-_{\ell}(n) e(n\bar{z}),
\end{align*}
where $ a^+_{\ell} (n) = 0$ for $n \ll 0$ and  $a^-_{\ell}(n) =0 $ for $n \gg 0$.
We let $H^+_0(\G)$ be the subspace of those $f$ that satisfy $ a^-_{\ell}(n) = 0$ for $n \geq 0$ (for all $\ell$). It consists for those $f \in H_0(\G)$ such that $f^-$ is exponentially decreasing at the cusps. It significance lies in the fact that $\xi_0$ maps $H_0^+(\G)$ onto $S_2(\G)$, the space of weight $2$ cusp forms for $\G$.
We define ${\bf{t}}_f(h,m)$ for $m \geq 0$ as before, while we define the modular trace of negative index ${\bf{t}}_f(h,-dm^2)$ by replacing $a_{\ell}(n)$ with the holomorphic coefficients $a^+_{\ell}(n)$.
\begin{theorem}
Let $ f \in H^+_0(\G)$ and assume that $a^+_{\ell}(0) =0$ for all $\ell$. Then
\[
I_h(\tau,f) = \sum_{m \geq 0} {\bf{t}}_f(h,m) q^m + \sum_{m > 0} {\bf{t}}_f(h,-dm^2) q^{-dm^2}.
\]
\end{theorem}

\begin{proof}

Since $f$ is harmonic, the proofs for the positive coefficients and for $q(X) = m \notin -d \left( \Q^{\times}\right)^2$, the non-split case, are still valid. That is, Propositions~\ref{prop1} and \ref{prop2a} carry over with no change.
  The term for $X=0$ stays also the same. Hence we only need to analyze the orbital integrals over  the isotropic lines and for the split case, $q(X) = -dm^2$.
For the extension of Proposition~\ref{prop2b}, we let $X \in L_{-dm^2}$ and see
\begin{align*}
\int_{M} & f(z) \sum_{\g \in \G} \varphi^0(\sqrt{v}X, \g z)
=
\frac1{2\pi i} \int_{M} f (z) \bar{\partial} \partial \sum_{\g \in \G} \xi^0(\sqrt{v}X, \g z) \\
&=
\frac1{2\pi i} \int_{M } d \left( f(z)  \partial \sum_{\g \in  \G} \xi^0(\sqrt{v}X, \g z) \right) - \frac1{2\pi i} \int_{M }   ( \bar{\partial}   f(z) ) \partial \sum_{\g \in  \G} \xi^0(\sqrt{v}X, \g z).
\end{align*}
The first term is handled in exactly the same manner as in the proof of Proposition~\ref{prop2b}. Only at the end of the proof of Lemma~\ref{britta}, when inserting the Fourier expansion of $f$, an extra term occurs. But one easily sees that this extra term vanishes in the limit. For the second term, we have
\begin{align*}
\int_{M }   ( \bar{\partial}   f(z) ) \partial \sum_{\g \in  \G} \xi^0(\sqrt{v}X, \g z) &=  -\int_{M } d \left( \bar{\partial} f(z) \sum_{\g \in  \G} \xi^0(\sqrt{v}X, \g z) \right) \\ &\quad + \int_{M }   ( \partial \bar{\partial}   f(z) ) \sum_{\g \in  \G} \xi^0(\sqrt{v}X, \g z).
\end{align*}
But the first summand vanishes by Stokes' theorem, since $  \bar{\partial} f(z)$ is rapidly decreasing as $f \in H^+_0$, while second term is zero since
$\partial \bar{\partial}   f(z) =0$ as $f$ is harmonic.
The orbital integrals over the isotropic lines are treated in the same manner.
\end{proof}

\end{document}